\documentclass[a4paper]{article}
\usepackage{amsfonts}
\usepackage{amsmath}
\usepackage{geometry}
\usepackage{graphicx}
\usepackage{multirow}

\newcommand*{\myalign}
[2]{\multicolumn{1}{#1}{#2}}
\input{tcilatex}
\begin{document}

\title{Analytical and numerical treatment of \\
the heat conduction equation obtained via \\
time-fractional distributed-order heat conduction law}
\author{Velibor \v{Z}eli%
\begin{footnote}
{Linn\'{e} FLOW Centre, KTH Mechanics, SE-100 44 Stockholm, Sweden, velibor@mech.kth.se},
\end{footnote}
Du\v{s}an Zorica%
\begin{footnote} {Mathematical Institute, Serbian
Academy of Arts and Sciences, Kneza Mihaila 36, 11000 Beograd,
Serbia, dusan\textunderscore zorica@mi.sanu.ac.rs and Department
of Physics, Faculty of Sciences, University of Novi Sad, Trg D.
Obradovi\'{c}a 3, 21000 Novi Sad, Serbia}
\end{footnote}}
\maketitle

\begin{abstract}
\noindent Generalization of the heat conduction equation is obtained by
considering the system of equations consisting of the energy balance
equation and fractional-order constitutive heat conduction law, assumed in
the form of the distributed-order Cattaneo type. The Cauchy problem for
system of energy balance equation and constitutive heat conduction law is
treated analytically through Fourier and Laplace integral transform methods,
as well as numerically by the method of finite differences through
Adams-Bashforth and Gr\"{u}nwald-Letnikov schemes for approximation
derivatives in temporal domain and leap frog scheme for spatial derivatives.
Numerical examples, showing time evolution of temperature and heat flux
spatial profiles, demonstrate applicability and good agreement of both
methods in cases of multi-term and power-type distributed-order heat
conduction laws.

\bigskip

\end{abstract}

\section{Introduction}

Heat conduction in one-dimensional rigid material is considered on infinite
spatial domain $x\in 
\mathbb{R}
$ and for time $t>0.$ Generalization of the heat conduction equation is
considered by treating two different processes in a material. The first one
is heating, described by the energy balance equation%
\begin{equation}
\rho c\frac{\partial }{\partial t}T\left( x,t\right) =-\frac{\partial }{%
\partial x}q\left( x,t\right) ,  \label{EB}
\end{equation}%
where $\rho $ is used to denote the material density, $c$ is the specific
heat capacity, while $T$ and $q$ denote temperature and heat flux
respectively. The second one is heat conduction, described by the Cattaneo
type time-fractional distributed-order heat conduction law%
\begin{equation}
\int_{0}^{1}\phi \left( \gamma \right) \,{}_{0}^{c}\mathrm{D}_{t}^{\gamma
}q\left( x,t\right) \mathrm{d}\gamma =-\lambda \frac{\partial }{\partial x}%
T\left( x,t\right) ,  \label{CE}
\end{equation}%
where $\phi $ is the constitutive function or distribution, $\lambda $ is
the thermal conductivity and ${}_{0}^{c}\mathrm{D}_{t}^{\gamma }$ is the
operator of Caputo fractional differentiation of order $\gamma \in \left(
0,1\right) ,$ defined by 
\begin{equation*}
{}_{0}^{c}\mathrm{D}_{t}^{\gamma }y\left( t\right) =\frac{t^{-\gamma }}{%
\Gamma \left( 1-\gamma \right) }\ast _{t}\frac{\mathrm{d}y\left( t\right) }{%
\mathrm{d}t},
\end{equation*}%
see \cite{TAFDE}, with $\ast _{t}$ denoting the convolution in time: $%
f\left( t\right) \ast _{t}g\left( t\right) =\int_{0}^{t}f\left( u\right)
g\left( t-u\right) \mathrm{d}u.$

Rather than obtaining and solving a single heat conduction equation, the aim
is to solve the system of equations consisting of energy balance equation (%
\ref{EB}) and constitutive equation (\ref{CE}), subject to initial 
\begin{equation}
T\left( x,0\right) =T_{0}\left( x\right) ,\;\;q\left( x,0\right) =0,\;\;x\in 
\mathbb{R}
,  \label{IC}
\end{equation}%
and boundary conditions%
\begin{equation}
\lim_{x\rightarrow \pm \infty }T\left( x,t\right) =0,\;\;\lim_{x\rightarrow
\pm \infty }q\left( x,t\right) =0,\;\;t>0.  \label{BC}
\end{equation}%
In particular, two special cases of the heat conduction law (\ref{CE}) will
be examined: multi-term heat conduction law, obtained for the choice of
constitutive distribution as%
\begin{equation}
\phi \left( \gamma \right) =\tau _{0}\,\delta \left( \gamma -\alpha
_{0}\right) +\sum_{\nu =1}^{N}\tau _{\nu }\,\delta \left( \gamma -\alpha
_{\nu }\right) ,\;\;%
\begin{array}{c}
0\leq \alpha _{0}<\ldots <\alpha _{N}<1, \\ 
\tau _{0},\tau _{1},\ldots ,\tau _{N}>0,%
\end{array}
\label{MT}
\end{equation}%
consisting of at least two terms, where $\delta $ is used to denote the
Dirac $\delta $-distribution, and power-type distributed-order heat
conduction law, obtained for the choice of constitutive function as%
\begin{equation}
\phi \left( \gamma \right) =\tau ^{\gamma },\;\;\tau >0.  \label{PTDO}
\end{equation}

The constitutive equation (\ref{CE}) represents the generalization of known
heat conduction laws such as Fourier, Cattaneo, fractional Cattaneo, which
are obtained by choosing the constitutive distribution as 
\begin{equation}
\phi \left( \gamma \right) =\delta \left( \gamma \right) ,\;\;\phi \left(
\gamma \right) =\tau \delta \left( \gamma -1\right) +\delta \left( \gamma
\right) ,\;\;\phi \left( \gamma \right) =\tau \delta \left( \gamma -\alpha
\right) +\delta \left( \gamma \right) ,  \label{CB}
\end{equation}%
respectively. The approach of considering generalized heat conduction
equation through system of balance and constitutive equation is also adopted
in \cite{Alvarez} within the classical theory using the analogy with
circuits and extending the results within the theory of fractional calculus
in \cite{FernandezAnaya}. Anomalous transport processes through space and
time fractional generalizations of the Cattaneo heat conduction law are
studied in \cite{CompteMetzler}. Time and space fractional heat conduction
of Cattaneo type is studied, analytically on infinite domain in \cite{AKOZ-1}
and with physical justification for non-locality introduction in \cite%
{BorinoDiPaolaZingales,MongioviZingales,Zingales}. Heat conduction problem
with the Riesz space fractional generalization of the Cattaneo-Christov heat
conduction model is numerically treated in \cite{LiuZhengLiuZhang}. Heat
conduction problems with different heat conduction laws in terms of the
classical theory are reviewed in \cite{Jos}, while in \cite{APSZ-1} there is
a collection of heat conduction problems within the theory of fractional
calculus.

By combining the energy balance equation (\ref{EB}) with the constitutive
equation (\ref{CE}), where constitutive distributions are given by (\ref{CB}%
), the classical heat conduction, telegraph and fractional telegraph
equations are obtained as 
\begin{equation}
\frac{\partial T}{\partial t}=\mathcal{D}\frac{\partial ^{2}T}{\partial x^{2}%
},\;\;\tau \frac{\partial ^{2}T}{\partial t^{2}}+\frac{\partial T}{\partial t%
}=\mathcal{D}\frac{\partial ^{2}T}{\partial x^{2}}\;\;\text{and}\;\;\tau
\,{}_{0}^{c}\mathrm{D}_{t}^{\alpha +1}T+\frac{\partial T}{\partial t}=%
\mathcal{D}\frac{\partial ^{2}T}{\partial x^{2}},  \label{HE}
\end{equation}%
respectively, with the thermal diffusivity $\mathcal{D=}\frac{\lambda }{\rho
c}.$ In \cite{Mohanty}, an unconditionally stable difference scheme is
developed for (\ref{HE})$_{2}$ having the source term, while the equation of
the form 
\begin{equation*}
\tau \,\frac{\partial ^{2}T}{\partial t^{2}}+{}_{0}^{c}\mathrm{D}%
_{t}^{\alpha }T=\mathcal{D}\frac{\partial ^{2}T}{\partial x^{2}},\;\;\text{%
with}\;\;\alpha \in \left( 1,2\right) ,
\end{equation*}%
is solved for the Dirichlet boundary conditions in \cite{QiGuo}. The
equation of the form (\ref{HE})$_{3}$, with the additional (source) term
corresponding to the impulse laser penetration into material causing its
heating, is solved analytically in \cite{QiXuGuo}, while in \cite{MishraRai}
the similar problem with time-exponential decay of laser heating was treated
numerically by developing unconditionally stable compact difference scheme.

Classical heat conduction equation (\ref{HE})$_{1}$ is often generalized
within the theory of fractional calculus by replacing the first order
partial time derivative with the fractional one, obtaining diffusion-wave
equation%
\begin{equation}
{}_{0}^{c}\mathrm{D}_{t}^{\alpha }T=\mathcal{D}\frac{\partial ^{2}T}{%
\partial x^{2}},  \label{DWE}
\end{equation}%
describing subdiffusion if $\alpha \in \left( 0,1\right) $ and
superdiffusion if $\alpha \in \left( 1,2\right) .$ The pioneering work on
diffusion-wave equation is \cite{Mai96}, followed by further explorations in 
\cite{GLM,MaLuPa}. Multi-dimensional variants of (\ref{DWE}) are studied in 
\cite{Hanyga-sf,Hanyga-stf,Hanyga-tf} for space-fractional,
space-time-fractional and time-fractional cases. Dirichlet problem for (\ref%
{DWE}), using several finite difference schemes, such as explicit, implicit
and Crank-Nicolson schemes is considered in \cite{ZecovaTerpak}. Similar
problem with the source term is numerically treated in \cite{SunWu}, while
an adaptive scheme was developed in \cite{YusteQuintanaMurillo}.
Diffusion-wave equation (\ref{DWE}) on bounded domain was numerically
considered in \cite{ZhuangLiu} through implicit difference scheme, while the
non-local variant of (\ref{DWE}) is analyzed in \cite{ShenLiuAnh} through
explicit and implicit finite difference schemes. In order to numerically
study anomalous infiltration in porous media the non-linear and
variable-order version of (\ref{DWE}) is used in \cite{ShenLiuLiuAnh}.

The fractional generalization of telegraph equation (\ref{HE})$_{2}$ can be
found in various forms different from (\ref{HE})$_{3}$ like 
\begin{equation}
\tau \,{}_{0}^{c}\mathrm{D}_{t}^{\alpha }T+{}_{0}^{c}\mathrm{D}_{t}^{\beta
}T=\mathcal{D}\frac{\partial ^{2}T}{\partial x^{2}},\;\;\text{or}\;\;\tau
\,{}_{0}^{c}\mathrm{D}_{t}^{2\alpha }T+{}_{0}^{c}\mathrm{D}_{t}^{\alpha }T=%
\mathcal{D}\frac{\partial ^{2}T}{\partial x^{2}},  \label{TEs}
\end{equation}%
with $0<\alpha <\beta <2,$ or $\alpha \in \left( 0,1\right) .$ Similarly as
in \cite{APZ}, where (\ref{TEs})$_{1}$ was treated analytically on infinite,
semi-infinite and finite domains, in \cite{XuQiJiang} the problem on the
semi-infinite domain for (\ref{TEs})$_{1}$ was treated for a special choice
of the boundary conditions. By the use of analytical methods, different
versions of telegraph equation are treated on unbounded and bounded domains
in \cite{CascavalEcksteinFrotaGoldstein,ChenLiuAnh,Huang}, including the
non-locality as in \cite{YakubovichRodrigues}. The non-local version of (\ref%
{TEs})$_{1}$ on unbounded domain, where the non-locality is expressed
through the fractional Laplacian is analytically treated in \cite%
{CamargoOliveiraVaz,CamargoCO08}, while the non-local version of the
telegraph equation (\ref{TEs})$_{2}$ on infinite domain is treated in \cite%
{QiJiang}.

Generalizing the fractional telegraph equation by adding terms containing
fractional derivatives of different orders lead to the distributed-order
diffusion-wave equation%
\begin{equation}
\int_{0}^{2}\phi \left( \gamma \right) \,{}_{0}^{c}\mathrm{D}_{t}^{\gamma
}T\,\mathrm{d}\gamma =\mathcal{D}\frac{\partial ^{2}T}{\partial x^{2}},
\label{DOE}
\end{equation}%
analytically analyzed in \cite%
{APZ-1,APZ-2,MaGoMi,MainardiPagniniGorenflo,MaPaMuGo}. A compact difference
scheme for (\ref{DOE}) on bounded domains, with source term included, is
developed in \cite{YeLiuAnh}, as well as in \cite{MorgadoRebelo}, where the
similar equation is also numerically analyzed. In \cite%
{MeerschaertNaneVellaisamy}, equation (\ref{DOE}) including the classical
Laplacian is considered on a bounded multi-dimensional domain. Local,
two-sided space-fractional, and Riesz space fractional variants of (\ref{DOE}%
) are analyzed through the implicit finite difference schemes in \cite%
{HuLiuAnhTurner,HuLiuTurnerAnh,YeLiuAnhTurnerIMA}. Multi-term
time-fractional diffusion type equation is considered in \cite{RJ}, while
the maximum principle and numerical method for the multi-term time-space
heat conduction equation of fractional order is considered in \cite%
{YeLiuAnhTurnerAMC}.

\section{Solution to Cauchy problem}

The Cauchy initial value problem on the real axis ($x\in 
\mathbb{R}
,$ $t>0$) for the time-fractional distributed-order Cattaneo type heat
conduction, i.e., the system of energy balance equation (\ref{EB}) and
constitutive Cattaneo type time-fractional distributed-order heat conduction
law (\ref{CE}), subject to initial (\ref{IC}) and boundary conditions (\ref%
{BC}), will be analytically solved by the means of integral transform
methods: Fourier transform with respect to spatial coordinate and Laplace
transform with respect to time, as well as by the finite difference method:
leap frog numerical scheme for spatial coordinate, along with Gr\"{u}%
nwald-Letnikov and third-order Adams-Bashforth temporal numerical schemes.
The use of Adams-Bashforth scheme will prove to give more accurate and
stable results when compared with centered, centered with RAW filter and
Euler schemes. Two cases of the constitutive equation (\ref{CE}) will be
examined: multi-term heat conduction law, with the constitutive distribution
given by (\ref{MT}), and power-type distributed-order heat conduction law,
with the constitutive function given by (\ref{PTDO}).

Dimensionless quantities%
\begin{equation*}
\bar{x}=\frac{x}{x^{\ast }},\;\;\bar{t}=\frac{t}{t^{\ast }},\;\;x^{\ast }=%
\sqrt{\frac{\lambda }{\rho c}t^{\ast }},\;\;\bar{T}=\frac{T}{\Theta _{0}}%
-1,\;\;\bar{q}=q\frac{1}{\Theta _{0}}\sqrt{\frac{t^{\ast }}{\lambda \rho c}}%
,\;\;\bar{\phi}=\frac{\phi }{\left( t^{\ast }\right) ^{\gamma }},
\end{equation*}%
where the time-scale $t^{\ast }$ will be determined according to the choice
of constitutive distribution/function $\phi ,$ see (\ref{t-star}) below, and
where the constant $\Theta _{0}$ represents the reference temperature,
introduced into system of equations (\ref{EB}) and (\ref{CE}), with
subsequent omittance of bars, yield the following form of governing
equations 
\begin{gather}
\frac{\partial }{\partial t}T\left( x,t\right) =-\frac{\partial }{\partial x}%
q\left( x,t\right) ,\;\;x\in 
\mathbb{R}
,\;t>0,  \label{EB-bd} \\
\int_{0}^{1}\phi \left( \gamma \right) \,{}_{0}^{c}\mathrm{D}_{t}^{\gamma
}q\left( x,t\right) \mathrm{d}\gamma =-\frac{\partial }{\partial x}T\left(
x,t\right) ,\;\;x\in 
\mathbb{R}
,\;t>0,  \label{CE-bd}
\end{gather}%
while the constitutive distribution (\ref{MT}) and the constitutive function
(\ref{PTDO}) become ($0\leq \alpha _{0}<\ldots <\alpha _{N}<1$)%
\begin{equation}
\phi \left( \gamma \right) =\delta \left( \gamma -\alpha _{0}\right)
+\sum_{\nu =1}^{N}\tau _{\nu }\,\delta \left( \gamma -\alpha _{\nu }\right)
\;\;\text{and}\;\;\phi \left( \gamma \right) =1,  \label{CDF-bd}
\end{equation}%
respectively. In the case of constitutive distribution (\ref{MT}),
respectively constitutive function (\ref{PTDO}), the time-scales%
\begin{equation}
t^{\ast }=\tau _{0}^{\frac{1}{\alpha _{0}}},\;\;\text{respectively}%
\;\;t^{\ast }=\tau ,  \label{t-star}
\end{equation}%
yield $\bar{\tau}_{\nu }=\tau _{\nu }\,\tau _{0}^{-\frac{\alpha _{\nu }}{%
\alpha _{0}}},$ $\nu =1,2,\ldots ,N,$ where bar is omitted in (\ref{CDF-bd}%
), respectively $\bar{\tau}=1.$

Governing equations (\ref{EB-bd}) and (\ref{CE-bd}), with constitutive
distribution/function (\ref{CDF-bd}), are subject to (dimensionless) initial
and boundary conditions 
\begin{gather}
T\left( x,0\right) =T_{0}\left( x\right) ,\;\;q\left( x,0\right) =0,\;\;x\in 
\mathbb{R}
,  \label{IC-bd} \\
\lim_{x\rightarrow \pm \infty }T\left( x,t\right) =0,\;\;\lim_{x\rightarrow
\pm \infty }q\left( x,t\right) =0,\;\;t>0.  \label{BC-bd}
\end{gather}

\subsection{Analytical solution}

Governing equations (\ref{EB-bd}) and (\ref{CE-bd}), with initial (\ref%
{IC-bd}) and boundary conditions (\ref{BC-bd}), will be analytically solved
by the integral transform method. Application of the Fourier, $\hat{f}\left(
\xi \right) =\mathcal{F}\left[ f\left( x\right) \right] \left( \xi \right)
=\int_{-\infty }^{\infty }f\left( x\right) \mathrm{e}^{-\mathrm{i}\xi x}%
\mathrm{d}x,$ $\xi \in 
\mathbb{R}
,$ and Laplace transform $\tilde{f}\left( s\right) =\mathcal{L}\left[
f\left( t\right) \right] \left( s\right) =\int_{0}^{\infty }f\left( t\right) 
\mathrm{e}^{-st}\mathrm{d}t,$ $\func{Re}s>0,$ to system of equations (\ref%
{EB-bd}) and (\ref{CE-bd}), with (\ref{IC-bd}) and (\ref{BC-bd}) taken into
account, yields%
\begin{gather}
s\widehat{\tilde{T}}\left( \xi ,s\right) -\hat{T}_{0}\left( \xi \right) =-%
\mathrm{i}\xi \,\widehat{\tilde{q}}\left( \xi ,s\right) ,\;\;\xi \in 
\mathbb{R}
,\;\func{Re}s>0,  \label{eb-lt} \\
\widehat{\tilde{q}}\left( \xi ,s\right) \Phi \left( s\right) =-\mathrm{i}\xi
\,\widehat{\tilde{T}}\left( \xi ,s\right) ,\;\;\xi \in 
\mathbb{R}
,\;\func{Re}s>0,  \label{ce-lt}
\end{gather}%
where 
\begin{equation*}
\Phi \left( s\right) =\int_{0}^{1}\phi \left( \gamma \right) s^{\gamma }%
\mathrm{d}\gamma ,\;\;\func{Re}s>0,
\end{equation*}%
in cases of constitutive distribution/function (\ref{CDF-bd}) takes the
following forms%
\begin{equation}
\Phi \left( s\right) =s^{\alpha _{0}}+\sum_{\nu =1}^{N}\tau _{\nu
}\,s^{\alpha _{\nu }}\;\;\text{and}\;\;\Phi \left( s\right) =\frac{s-1}{\ln s%
},\;\;\func{Re}s>0.  \label{cdf-lt}
\end{equation}

Solution to system of equations (\ref{eb-lt}) and (\ref{ce-lt}) with respect
to $\widehat{\tilde{T}}$ and $\widehat{\tilde{q}}$ reads ($\xi \in 
\mathbb{R}
,$\ $\func{Re}s>0$)%
\begin{equation}
\widehat{\tilde{T}}\left( \xi ,s\right) =\hat{T}_{0}\left( \xi \right) \frac{%
\Phi \left( s\right) }{\xi ^{2}+s\Phi \left( s\right) }\;\;\text{and}\;\;%
\widehat{\tilde{q}}\left( \xi ,s\right) =-\hat{T}_{0}\left( \xi \right) 
\frac{\mathrm{i}\xi }{\xi ^{2}+s\Phi \left( s\right) }.
\label{T,q-tilda-het}
\end{equation}%
Using the well-known Fourier inversion formula%
\begin{equation}
\mathcal{F}^{-1}\left[ \frac{1}{\xi ^{2}+\lambda }\right] \left( x\right) =%
\frac{1}{2\sqrt{\lambda }}\mathrm{e}^{-\left\vert x\right\vert \sqrt{\lambda 
}},\;\;x\in 
\mathbb{R}
,\;\lambda \in 
\mathbb{C}
\backslash \left( -\infty ,0\right] ,  \label{FIF}
\end{equation}%
in (\ref{T,q-tilda-het}), along with the Fourier transform of a derivative
and convolution, one obtains%
\begin{equation}
T\left( x,t\right) =T_{0}\left( x\right) \ast _{x}P\left( x,t\right) \;\;%
\text{and}\;\;q\left( x,t\right) =T_{0}\left( x\right) \ast _{x}Q\left(
x,t\right) ,  \label{T,q}
\end{equation}%
after additional inversion of the Laplace transform, where ($x\in 
\mathbb{R}
,$\ $\func{Re}s>0$)%
\begin{eqnarray}
\tilde{P}\left( x,s\right) &=&\frac{1}{2}\sqrt{\frac{\Phi \left( s\right) }{s%
}}\mathrm{e}^{-\left\vert x\right\vert \sqrt{s\Phi \left( s\right) }},
\label{P-tilde} \\
\tilde{Q}\left( x,s\right) &=&-\frac{1}{2\sqrt{s\Phi \left( s\right) }}\frac{%
\mathrm{d}}{\mathrm{d}x}\mathrm{e}^{-\left\vert x\right\vert \sqrt{s\Phi
\left( s\right) }}=\frac{1}{2}\mathrm{e}^{-\left\vert x\right\vert \sqrt{%
s\Phi \left( s\right) }}\limfunc{sgn}x,  \label{Q-tilde}
\end{eqnarray}%
with $\limfunc{sgn}x=2H\left( x\right) -1,$ $x\in 
\mathbb{R}
,$ and $H$ being the Heaviside function. The justification for using the
Fourier inversion formula (\ref{FIF}), as well as the argumentation that
complex square root is well-defined, is given in Appendix \ref{app}.

When the Laplace inversion formula 
\begin{equation*}
f\left( t\right) =\frac{1}{2\pi \mathrm{i}}\int_{c-\mathrm{i}\infty }^{c+%
\mathrm{i}\infty }\tilde{f}\left( s\right) \mathrm{e}^{st}\mathrm{d}%
s,\;\;t>0,
\end{equation*}%
is applied to $\tilde{P}$ and $\tilde{Q},$ given by (\ref{P-tilde}) and (\ref%
{Q-tilde}), then solution kernels $P$ and $Q,$ appearing in (\ref{T,q}), are
obtained for $x\in 
\mathbb{R}
,$ $t>0,$ as 
\begin{eqnarray}
P\left( x,t\right) &=&\frac{1}{4\pi }\int_{0}^{\infty }\left( \sqrt{\Phi
^{+}\left( p\right) }\,\mathrm{e}^{-\mathrm{i}\left\vert x\right\vert \sqrt{%
p\Phi ^{+}\left( p\right) }}+\sqrt{\Phi ^{-}\left( p\right) }\,\mathrm{e}^{%
\mathrm{i}\left\vert x\right\vert \sqrt{p\Phi ^{-}\left( p\right) }}\right) 
\frac{\mathrm{e}^{-pt}}{\sqrt{p}}\mathrm{d}p,  \label{P} \\
Q\left( x,t\right) &=&\frac{\func{sgn}x}{4\pi \mathrm{i}}\int_{0}^{\infty
}\left( \mathrm{e}^{\mathrm{i}\left\vert x\right\vert \sqrt{p\Phi ^{-}\left(
p\right) }}-\mathrm{e}^{-\mathrm{i}\left\vert x\right\vert \sqrt{p\Phi
^{+}\left( p\right) }}\right) \mathrm{e}^{-pt}\mathrm{d}p,  \label{Q}
\end{eqnarray}%
where 
\begin{equation*}
\Phi ^{+}\left( p\right) =\Phi \left( p\,\mathrm{e}^{\mathrm{i}\pi }\right)
\;\;\text{and}\;\;\Phi ^{-}\left( p\right) =\Phi \left( p\,\mathrm{e}^{-%
\mathrm{i}\pi }\right) ,
\end{equation*}%
with $\Phi $ given by (\ref{cdf-lt}).

\begin{figure}[tbph]
\centering
\includegraphics[scale=0.55]{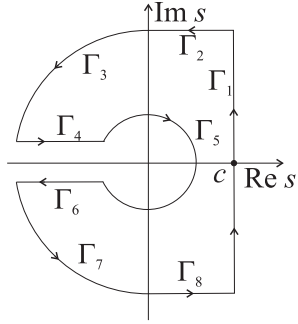}
\caption{Contour $\Gamma $.}
\label{hk}
\end{figure}

In order to obtain functions $P$ and $Q,$ (\ref{P}) and (\ref{Q}), the
Cauchy integral theorem 
\begin{equation}
\oint_{\Gamma }\tilde{P}\left( x,s\right) \mathrm{e}^{st}\mathrm{d}s=0,
\label{kit-p}
\end{equation}%
is used, where $\Gamma =\Gamma _{1}\cup \Gamma _{2}\cup \Gamma _{3}\cup
\Gamma _{4}\cup \Gamma _{5}\cup \Gamma _{6}\cup \Gamma _{7}\cup \Gamma _{8}$
is the closed contour shown in Figure \ref{hk}, within $\tilde{P}\left(
x,s\right) \mathrm{e}^{st},$ $x\in 
\mathbb{R}
,$ $t>0,$ is an analytic function, since apart from $s=0,$ there are no
other branching points of $\tilde{P}$ and $\tilde{Q}.$ 
The contour integration will be performed for $\tilde{P}$ only, since the
similar procedure and arguments are applicable for $\tilde{Q}$ as well. This
will be shown in the sequel by proving that $\Phi ,$ given by (\ref{cdf-lt}%
), has no zeros in the principal branch, i.e., for $\arg s\in \left( -\pi
,\pi \right) .$

The real and imaginary parts of $\Phi ,$ given by (\ref{cdf-lt})$_{1}$,
after substitution $s=\rho \,\mathrm{e}^{\mathrm{i}\varphi },$ $\rho >0,$ $%
\varphi \in \left( -\pi ,\pi \right) ,$ read%
\begin{eqnarray}
\func{Re}\Phi \left( \rho ,\varphi \right) &=&\rho ^{\alpha _{0}}\cos \left(
\alpha _{0}\varphi \right) +\sum_{\nu =1}^{N}\tau _{\nu }\rho ^{\alpha _{\nu
}}\cos \left( \alpha _{\nu }\varphi \right) ,  \label{fi-re} \\
\func{Im}\Phi \left( \rho ,\varphi \right) &=&\rho ^{\alpha _{0}}\sin \left(
\alpha _{0}\varphi \right) +\sum_{\nu =1}^{N}\tau _{\nu }\rho ^{\alpha _{\nu
}}\sin \left( \alpha _{\nu }\varphi \right) .  \label{fi-im}
\end{eqnarray}%
Since, by (\ref{fi-im}), it holds that $\func{Im}\Phi \left( \bar{s}\right)
=-\func{Im}\Phi \left( s\right) ,$ where bar denotes the complex
conjugation, it is sufficient to analyze function $\Phi $ for $\varphi \in %
\left[ 0,\pi \right) $ only. If $\varphi \in \left( 0,\pi \right) ,$ then $%
\func{Im}\Phi \left( \rho ,\varphi \right) >0,$ since for $0\leq \alpha
_{0}<\ldots <\alpha _{N}<1,$ it is valid that $\sin \left( \alpha _{\nu
}\varphi \right) >0,$ $\nu =0,1,\ldots ,N$ (except for possible $\alpha
_{0}=0,$ which does affect that $\func{Im}\Phi \left( \rho ,\varphi \right)
>0$). Note that $\func{Im}\Phi \left( \rho ,\varphi \right) =0,$ only if
there is a single term with $\alpha _{0}=0$ in the constitutive distribution
(\ref{CDF-bd})$_{1}$. This, however is not the case. If $\varphi =0,$ then,
by (\ref{fi-re}), $\func{Re}\Phi \left( \rho ,\varphi \right) >0.$
Therefore, function $\Phi ,$ (\ref{cdf-lt})$_{1}$, has no zeros in the
principal branch.

In the case of $\Phi ,$ given by (\ref{cdf-lt})$_{2},$ the corresponding
real and imaginary parts are%
\begin{eqnarray}
\func{Re}\Phi \left( \rho ,\varphi \right) &=&\frac{\ln \rho \left( \rho
\cos \varphi -1\right) +\rho \varphi \sin \varphi }{\ln ^{2}\rho +\varphi
^{2}},  \label{psi-re-1} \\
\func{Im}\Phi \left( \rho ,\varphi \right) &=&\frac{\rho \ln \rho \sin
\varphi -\varphi \left( \rho \cos \varphi -1\right) }{\ln ^{2}\rho +\varphi
^{2}}.  \label{psi-im-1}
\end{eqnarray}%
In order to determine zeros of $\Phi $, both real and imaginary part of $%
\Phi $ have to be zero, yielding the system of equations%
\begin{gather*}
\ln \rho \left( \rho \cos \varphi -1\right) +\rho \varphi \sin \varphi =0, \\
\rho \ln \rho \sin \varphi -\varphi \left( \rho \cos \varphi -1\right) =0.
\end{gather*}%
Combining the equations, one obtains%
\begin{equation*}
\left( \ln ^{2}\rho +\varphi ^{2}\right) \left( \rho \cos \varphi -1\right)
=0.
\end{equation*}%
If $\ln ^{2}\rho +\varphi ^{2}=0,$ then real and imaginary part of $\Phi $
are not well defined, see (\ref{psi-re-1}) and (\ref{psi-im-1}). Therefore, $%
\rho \cos \varphi -1=0,$ implying $\rho \varphi \sin \varphi =0,$ whose
solutions are $\rho =0,$ $\varphi =0,$ $\varphi =\pi .$ The first and third
solution are in contradiction with $\rho \cos \varphi -1=0$, while $\rho =1$
corresponds to the second one. By substituting $s=\left. \rho \,\mathrm{e}^{%
\mathrm{i}\varphi }\right\vert _{\varphi =0}$ in (\ref{cdf-lt})$_{2}$ one
obtains $\frac{\rho -1}{\ln \rho }=0,$ while, on the other hand, $\lim_{\rho
\rightarrow 1}\frac{\rho -1}{\ln \rho }=1,$ implying that $\rho =1$ is not a
zero of $\Phi .$ Thus, there are no zeros in the principal branch of
function $\Phi ,$ (\ref{cdf-lt})$_{2}$, as well.

The integrals along contours $\Gamma _{1},$ $\Gamma _{4},$ parametrized by $%
s=p\,\mathrm{e}^{\mathrm{i}\pi },$ and $\Gamma _{6},$ parametrized by $s=p\,%
\mathrm{e}^{-\mathrm{i}\pi },$ in the limit when $R$ tends to infinity and $%
r $ tends to zero, yield ($x\in 
\mathbb{R}
,$ $t>0$)%
\begin{eqnarray*}
\lim_{R\rightarrow \infty }\int_{\Gamma _{1}}\tilde{P}\left( x,s\right) 
\mathrm{e}^{st}\mathrm{d}s &=&2\pi \mathrm{i}P\left( x,t\right) , \\
\lim_{\substack{ R\rightarrow \infty ,  \\ r\rightarrow 0}}\int_{\Gamma _{4}}%
\tilde{P}\left( x,s\right) \mathrm{e}^{st}\mathrm{d}s &=&-\frac{1}{2\mathrm{i%
}}\lim_{\substack{ R\rightarrow \infty ,  \\ r\rightarrow 0}}\int_{R}^{r}%
\frac{\sqrt{\Phi \left( p\,\mathrm{e}^{\mathrm{i}\pi }\right) }}{\sqrt{p}}%
\mathrm{e}^{-\mathrm{i}\left\vert x\right\vert \sqrt{p\Phi \left( p\,\mathrm{%
e}^{\mathrm{i}\pi }\right) }}\mathrm{e}^{-pt}\mathrm{d}p \\
&=&\frac{1}{2\mathrm{i}}\int_{0}^{\infty }\frac{\sqrt{\Phi ^{+}\left(
p\right) }}{\sqrt{p}}\mathrm{e}^{-\mathrm{i}\left\vert x\right\vert \sqrt{%
p\Phi ^{+}\left( p\right) }}\mathrm{e}^{-pt}\mathrm{d}p, \\
\lim_{\substack{ R\rightarrow \infty ,  \\ r\rightarrow 0}}\int_{\Gamma _{6}}%
\tilde{P}\left( x,s\right) \mathrm{e}^{st}\mathrm{d}s &=&\frac{1}{2\mathrm{i}%
}\lim_{\substack{ R\rightarrow \infty ,  \\ r\rightarrow 0}}\int_{r}^{R}%
\frac{\sqrt{\Phi \left( p\,\mathrm{e}^{-\mathrm{i}\pi }\right) }}{\sqrt{p}}%
\mathrm{e}^{\mathrm{i}\left\vert x\right\vert \sqrt{p\Phi \left( p\,\mathrm{e%
}^{-\mathrm{i}\pi }\right) }}\mathrm{e}^{-pt}\mathrm{d}p \\
&=&\frac{1}{2\mathrm{i}}\int_{0}^{\infty }\frac{\sqrt{\Phi ^{-}\left(
p\right) }}{\sqrt{p}}\mathrm{e}^{\mathrm{i}\left\vert x\right\vert \sqrt{%
p\Phi ^{-}\left( p\right) }}\mathrm{e}^{-pt}\mathrm{d}p,
\end{eqnarray*}%
while the integrals along $\Gamma _{2},$ $\Gamma _{3},$ $\Gamma _{5},$ $%
\Gamma _{7}$ and $\Gamma _{8}$ are zero. Using aforementioned integrals in
the Cauchy integral theorem (\ref{kit-p}) yields solution kernel $P$ in the
form given by (\ref{P}).

The absolute value of an integral along contour $\Gamma _{2},$ parametrized
by $s=q+\mathrm{i}R,$ with $q\in \left( 0,c\right) $ and $R$ tending to
infinity, is estimated as%
\begin{equation}
\left\vert \int_{\Gamma _{2}}\tilde{P}\left( x,s\right) \mathrm{e}^{st}%
\mathrm{d}s\right\vert \leq \frac{1}{2}\int_{0}^{c}\frac{\left\vert \sqrt{%
\Phi \left( q+\mathrm{i}R\right) }\right\vert }{\left\vert \sqrt{q+\mathrm{i}%
R}\right\vert }\left\vert \mathrm{e}^{-\left\vert x\right\vert \sqrt{q+%
\mathrm{i}R}\sqrt{\Phi \left( q+\mathrm{i}R\right) }}\right\vert \mathrm{e}%
^{qt}\mathrm{d}q.  \label{int-gama-2}
\end{equation}%
Note that $\sqrt{q+\mathrm{i}R}\sim \sqrt{R}\mathrm{e}^{\mathrm{i}\frac{\pi 
}{4}},$ for $R\rightarrow \infty .$ For $\Phi $ given by (\ref{cdf-lt})$_{1}$%
, according to (\ref{fi-re}) and (\ref{fi-im}) and having in mind $0\leq
\alpha _{0}<\ldots <\alpha _{N}<1,$ as $R\rightarrow \infty ,$ one has%
\begin{gather*}
\func{Re}\Phi \left( R,\frac{\pi }{2}\right) \sim R^{\alpha _{N}}\cos \frac{%
\alpha _{N}\pi }{2}\;\;\text{and}\;\;\func{Im}\Phi \left( R,\frac{\pi }{2}%
\right) \sim R^{\alpha _{N}}\sin \frac{\alpha _{N}\pi }{2},\;\;\text{i.e.,}
\\
\sqrt{\Phi \left( q+\mathrm{i}R\right) }\sim R^{\frac{\alpha _{N}}{2}}%
\mathrm{e}^{\mathrm{i}\frac{\alpha _{N}\pi }{4}},
\end{gather*}%
while for $\Phi $ given by (\ref{cdf-lt})$_{2}$, according to (\ref{psi-re-1}%
) and (\ref{psi-im-1}), as $R\rightarrow \infty ,$ one has%
\begin{gather*}
\func{Re}\Phi \left( R,\frac{\pi }{2}\right) \sim \frac{\pi }{2}\frac{R}{\ln
^{2}R}\;\;\text{and}\;\;\func{Im}\Phi \left( R,\frac{\pi }{2}\right) \sim 
\frac{R}{\ln R},\;\;\text{i.e.,} \\
\sqrt{\Phi \left( q+\mathrm{i}R\right) }\sim \sqrt{\frac{R}{\ln R}}\mathrm{e}%
^{\mathrm{i}\frac{\pi }{4}},
\end{gather*}%
so that, as $R\rightarrow \infty ,$ in the first case of $\Phi ,$ (\ref%
{int-gama-2}) becomes%
\begin{equation*}
\left\vert \int_{\Gamma _{2}}\tilde{P}\left( x,s\right) \mathrm{e}^{st}%
\mathrm{d}s\right\vert \leq \frac{1}{2}\int_{0}^{c}\frac{1}{R^{\frac{%
1-\alpha _{N}}{2}}}\mathrm{e}^{-\left\vert x\right\vert R^{\frac{1+\alpha
_{N}}{2}}\cos \frac{\left( 1+\alpha _{N}\right) \pi }{4}}\mathrm{e}^{qt}%
\mathrm{d}q\leq 0,
\end{equation*}%
since $\cos \frac{\left( 1+\alpha _{N}\right) \pi }{4}>0,$ while in the
second case of $\Phi ,$ (\ref{int-gama-2}) becomes%
\begin{equation*}
\left\vert \int_{\Gamma _{2}}\tilde{P}\left( x,s\right) \mathrm{e}^{st}%
\mathrm{d}s\right\vert \leq \frac{1}{2}\int_{0}^{c}\frac{1}{\sqrt{\ln R}}%
\mathrm{e}^{qt}\mathrm{d}q\leq 0,
\end{equation*}%
implying in both cases $\lim_{R\rightarrow \infty }\int_{\Gamma _{2}}\tilde{P%
}\left( x,s\right) \mathrm{e}^{st}\mathrm{d}s=0.$ Similar arguments yield $%
\lim_{R\rightarrow \infty }\int_{\Gamma _{8}}\tilde{P}\left( x,s\right) 
\mathrm{e}^{st}\mathrm{d}s=0.$

On the contour $\Gamma _{3},$ parametrized by $s=R\,\mathrm{e}^{\mathrm{i}%
\varphi }$, $\varphi \in \left( \frac{\pi }{2},\pi \right) ,$ with $%
R\rightarrow \infty ,$ the absolute value of the corresponding integral is
estimated as%
\begin{equation}
\left\vert \int_{\Gamma _{3}}\tilde{P}\left( x,s\right) \mathrm{e}^{st}%
\mathrm{d}s\right\vert \leq \frac{1}{2}\int_{\frac{\pi }{2}}^{\pi }R\frac{%
\left\vert \sqrt{\Phi \left( R,\varphi \right) }\right\vert }{\sqrt{R}}%
\left\vert \mathrm{e}^{-\left\vert x\right\vert \sqrt{R\,\mathrm{e}^{\mathrm{%
i}\varphi }}\sqrt{\Phi \left( R,\varphi \right) }}\right\vert \mathrm{e}%
^{Rt\cos \varphi }\mathrm{d}\varphi .  \label{int-gama-3}
\end{equation}%
For $\Phi $ given by (\ref{cdf-lt})$_{1}$, according to (\ref{fi-re}) and (%
\ref{fi-im}), as $R\rightarrow \infty ,$ one has%
\begin{gather*}
\func{Re}\Phi \left( R,\varphi \right) \sim R^{\alpha _{N}}\cos \left(
\alpha _{N}\varphi \right) \;\;\text{and}\;\;\func{Im}\Phi \left( R,\varphi
\right) \sim R^{\alpha _{N}}\sin \left( \alpha _{N}\varphi \right) ,\;\;%
\text{i.e.,} \\
\sqrt{\Phi \left( R,\varphi \right) }\sim R^{\frac{\alpha _{N}}{2}}\mathrm{e}%
^{\mathrm{i}\frac{\alpha _{N}\varphi }{2}},
\end{gather*}%
while for $\Phi $ given by (\ref{cdf-lt})$_{2}$, according to (\ref{psi-re-1}%
) and (\ref{psi-im-1}), as $R\rightarrow \infty ,$ one has 
\begin{gather*}
\func{Re}\Phi \left( R,\varphi \right) \sim \frac{R}{\ln R}\cos \varphi \;\;%
\text{and}\;\;\func{Im}\Phi \left( R,\varphi \right) \sim \frac{R}{\ln R}%
\sin \varphi ,\;\;\text{i.e.,} \\
\sqrt{\Phi \left( R,\varphi \right) }\sim \sqrt{\frac{R}{\ln R}}\mathrm{e}^{%
\mathrm{i}\frac{\varphi }{2}},
\end{gather*}%
so that, as $R\rightarrow \infty ,$ in the first case of $\Phi ,$ (\ref%
{int-gama-3}) becomes%
\begin{equation*}
\left\vert \int_{\Gamma _{3}}\tilde{P}\left( x,s\right) \mathrm{e}^{st}%
\mathrm{d}s\right\vert \leq \frac{1}{2}\int_{\frac{\pi }{2}}^{\pi }R^{\frac{%
1+\alpha _{N}}{2}}\mathrm{e}^{Rt\cos \varphi -\left\vert x\right\vert R^{%
\frac{1+\alpha _{N}}{2}}\cos \frac{\left( 1+\alpha _{N}\right) \varphi }{2}}%
\mathrm{d}\varphi \leq 0,
\end{equation*}%
while in the second case of $\Phi ,$ (\ref{int-gama-3}) becomes%
\begin{equation*}
\left\vert \int_{\Gamma _{3}}\tilde{P}\left( x,s\right) \mathrm{e}^{st}%
\mathrm{d}s\right\vert \leq \frac{1}{2}\int_{\frac{\pi }{2}}^{\pi }\frac{R}{%
\sqrt{\ln R}}\mathrm{e}^{Rt\cos \varphi -\left\vert x\right\vert \frac{R}{%
\sqrt{\ln R}}\cos \varphi }\mathrm{d}\varphi \leq 0.
\end{equation*}%
In both cases $\cos \varphi <0$ and the first term in the exponential is of
the highest order, implying $\lim_{R\rightarrow \infty }\int_{\Gamma _{3}}%
\tilde{P}\left( x,s\right) \mathrm{e}^{st}\mathrm{d}s=0.$ Similar arguments
yield $\lim_{R\rightarrow \infty }\int_{\Gamma _{7}}\tilde{P}\left(
x,s\right) \mathrm{e}^{st}\mathrm{d}s=0.$

The absolute value of the integral along the contour $\Gamma _{5},$
parametrized by $s=r\,\mathrm{e}^{\mathrm{i}\varphi }$, $\varphi \in \left(
-\pi ,\pi \right) ,$ with $r\rightarrow 0,$ is estimated as%
\begin{equation}
\left\vert \int_{\Gamma _{5}}\tilde{P}\left( x,s\right) \mathrm{e}^{st}%
\mathrm{d}s\right\vert \leq \frac{1}{2}\int_{-\pi }^{\pi }r\frac{\left\vert 
\sqrt{\Phi \left( r,\varphi \right) }\right\vert }{\sqrt{r}}\left\vert 
\mathrm{e}^{-\left\vert x\right\vert \sqrt{r\,\mathrm{e}^{\mathrm{i}\varphi }%
}\sqrt{\Phi \left( r,\varphi \right) }}\right\vert \mathrm{e}^{rt\cos
\varphi }\mathrm{d}\varphi .  \label{int-gama-5}
\end{equation}%
For $\Phi $ given by (\ref{cdf-lt})$_{1}$, according to (\ref{fi-re}) and (%
\ref{fi-im}), as $r\rightarrow 0,$ one has%
\begin{gather*}
\func{Re}\Phi \left( r,\varphi \right) \sim r^{\alpha _{0}}\cos \left(
\alpha _{0}\varphi \right) \;\;\text{and}\;\;\func{Im}\Phi \left( r,\varphi
\right) \sim r^{\alpha _{0}}\sin \left( \alpha _{0}\varphi \right) ,\;\;%
\text{i.e.,} \\
\sqrt{\Phi \left( r,\varphi \right) }\sim r^{\frac{\alpha _{0}}{2}}\mathrm{e}%
^{\mathrm{i}\frac{\alpha _{0}\varphi }{2}},
\end{gather*}%
while for $\Phi $ given by (\ref{cdf-lt})$_{2}$, according to (\ref{psi-re-1}%
) and (\ref{psi-im-1}), as $r\rightarrow 0,$ one has%
\begin{gather*}
\func{Re}\Phi \left( r,\varphi \right) \sim -\frac{1}{\ln r}\;\;\text{and}%
\;\;\func{Im}\Phi \left( r,\varphi \right) \sim \frac{\varphi }{\ln ^{2}r}%
,\;\;\text{i.e.,} \\
\sqrt{\Phi \left( r,\varphi \right) }\sim \frac{1}{\sqrt{\left\vert \ln
r\right\vert }},
\end{gather*}%
so that, as $r\rightarrow 0,$ in the first case of $\Phi ,$ (\ref{int-gama-5}%
) becomes%
\begin{equation*}
\left\vert \int_{\Gamma _{5}}\tilde{P}\left( x,s\right) \mathrm{e}^{st}%
\mathrm{d}s\right\vert \leq \frac{1}{2}\int_{-\pi }^{\pi }r^{\frac{1+\alpha
_{0}}{2}}\mathrm{e}^{rt\cos \varphi -\left\vert x\right\vert r^{\frac{%
1+\alpha _{0}}{2}}\cos \frac{\left( 1+\alpha _{0}\right) \varphi }{2}}%
\mathrm{d}\varphi \leq 0,
\end{equation*}%
while in the second case of $\Phi ,$ (\ref{int-gama-5}) becomes%
\begin{equation*}
\left\vert \int_{\Gamma _{5}}\tilde{P}\left( x,s\right) \mathrm{e}^{st}%
\mathrm{d}s\right\vert \leq \frac{1}{2}\int_{-\pi }^{\pi }\sqrt{\frac{r}{%
\left\vert \ln r\right\vert }}\mathrm{e}^{rt\cos \varphi -\left\vert
x\right\vert \sqrt{\frac{r}{\left\vert \ln r\right\vert }}\cos \frac{\varphi 
}{2}}\mathrm{d}\varphi \leq 0,
\end{equation*}%
implying in both cases $\lim_{r\rightarrow 0}\int_{\Gamma _{5}}\tilde{P}%
\left( x,s\right) \mathrm{e}^{st}\mathrm{d}s=0.$

\subsection{Numerical solution through finite difference scheme \label{fds}}

The Cauchy problem consisting of the (dimensionless) energy balance equation
(\ref{EB-bd}) and (dimensionless) constitutive Cattaneo type time-fractional
distributed-order heat conduction law (\ref{CE-bd}), corresponding to the
time-fractional distributed-order Cattaneo type heat conduction, subject to
initial (\ref{IC-bd}) and boundary conditions (\ref{BC-bd}), will be
numerically solved, as the coupled system of equations, through the finite
difference method, where discretization takes place on the spatial and
temporal domain, with $\Delta x$ and $\Delta t$ being their steps,
respectively.

The spatial derivatives will be approximated using the leap frog scheme%
\begin{equation}
\left. \frac{\partial }{\partial x}y\left( x,t\right) \right\vert
_{x=j\Delta x,\;t=n\Delta t}\approx \frac{y_{j+1}^{n}-y_{j-1}^{n}}{2\Delta x}%
,  \label{lf}
\end{equation}%
which is standardly used second order accuracy scheme. The numerical
calculations have shown that, unless using this scheme in discretization of
all spatial derivatives appearing in the governing equations (\ref{EB-bd})
and (\ref{CE-bd}), the divergence occurs during the calculation of the
solution. The possibility of using other schemes in discretization of all
spatial derivatives in the governing equations was not considered.

The third order accuracy Adams-Bashforth scheme for the first order
differential equation of the type%
\begin{equation*}
\frac{\partial }{\partial t}y\left( x,t\right) =f\left( x,t\right) ,
\end{equation*}%
reads%
\begin{equation}
\frac{y_{j}^{n+1}-y_{j}^{n}}{\Delta t}=\frac{1}{12}\left(
23f_{j}^{n}-16f_{j}^{n-1}+5f_{j}^{n-2}\right) ,  \label{ab}
\end{equation}%
and it will be used to approximate the energy balance equation (\ref{EB-bd}%
), since it significantly dumps the computational mode for $\Delta t$ small
enough, see \cite{Durran}. The use of the third order Adams-Bashforth scheme
will be justified in Section \ref{cds} by comparing its performance in
stability and accuracy while calculating the solution to governing
equations, with the first order accuracy Euler scheme%
\begin{equation}
\frac{y_{j}^{n+1}-y_{j}^{n}}{\Delta t}=f_{j}^{n},  \label{es}
\end{equation}%
the second order accuracy centered scheme 
\begin{equation}
\frac{y_{j}^{n+1}-y_{j}^{n-1}}{2\Delta t}=f_{j}^{n},  \label{cs}
\end{equation}%
and centered scheme with RAW filter, which is of the third order accuracy,
thus improved when compared with the centered scheme. For the properties and
implementation of the centered scheme with RAW filter see \cite{Williams}.

Although the Caputo fractional derivative appears in the constitutive
equation (\ref{CE-bd}), Gr\"{u}nwald-Letnikov approximation of the
Riemann-Liouville derivative will be used, due to the zero initial condition
(\ref{IC-bd})$_{2}$, implying the equivalence of the Caputo and
Riemann-Liouville fractional derivatives, so that, for $\alpha \in \left(
0,1\right) $ 
\begin{equation*}
\left[ {}_{0}^{c}\mathrm{D}_{t}^{\gamma }y\left( x,t\right) \right]
_{x=j\Delta x,\;t=n\Delta t}\approx \frac{1}{\left( \Delta t\right) ^{\gamma
}}y_{j}^{n}+\frac{1}{\left( \Delta t\right) ^{\gamma }}\sum_{k=1}^{n}\omega
_{k}\left( \gamma \right) y_{j}^{n-k},\;\;\omega _{k}\left( \gamma \right)
=\left( -1\right) ^{k}\binom{\gamma }{k},
\end{equation*}%
see \cite{pod}. When applying the Gr\"{u}nwald-Letnikov approximation, in
order to avoid calculation of the binomial coefficients, i.e., gamma
functions for large arguments, recurrence relation 
\begin{equation*}
\omega _{k}\left( \gamma \right) =\left( 1-\frac{1+\gamma }{k}\right) \omega
_{k-1}\left( \gamma \right) ,\;\;\omega _{0}\left( \gamma \right) =1,
\end{equation*}%
is adopted. Additionally, the integral appearing in the distributed-order
constitutive law (\ref{CE-bd}) will be approximated using the trapezoidal
method.

Employing the leap frog scheme (\ref{lf}) in spatial domain and
Adams-Bashforth scheme (\ref{ab}) in the energy balance equation (\ref{EB-bd}%
), the approximation of governing equations (\ref{EB-bd}), (\ref{CE-bd})
reads%
\begin{gather*}
\frac{T_{j}^{n+1}-T_{j}^{n}}{\Delta t}=-\frac{1}{12}\left( 23\frac{%
q_{j+1}^{n}-q_{j-1}^{n}}{2\Delta x}-16\frac{q_{j+1}^{n-1}-q_{j-1}^{n-1}}{%
2\Delta x}+5\frac{q_{j+1}^{n-2}-q_{j-1}^{n-2}}{2\Delta x}\right) , \\
\int_{0}^{1}\phi \left( \gamma \right) \left( \frac{1}{\left( \Delta
t\right) ^{\gamma }}q_{j}^{n}+\frac{1}{\left( \Delta t\right) ^{\gamma }}%
\sum_{k=1}^{n}\omega _{k}\left( \gamma \right) q_{j}^{n-k}\right) \mathrm{d}%
\gamma =-\frac{T_{j+1}^{n}-T_{j-1}^{n}}{2\Delta x},
\end{gather*}%
with the initial conditions (\ref{IC-bd}) yielding%
\begin{equation}
T_{j}^{0}=\left( T_{0}\right) _{j},\;\;q_{j}^{0}=0,  \label{ic-a}
\end{equation}%
so that one obtains%
\begin{eqnarray}
T_{j}^{n+1} &=&T_{j}^{n}-\frac{\Delta t}{24\Delta x}  \notag \\
&&\times \left( 23\left( q_{j+1}^{n}-q_{j-1}^{n}\right) -16\left(
q_{j+1}^{n-1}-q_{j-1}^{n-1}\right) +5\left(
q_{j+1}^{n-2}-q_{j-1}^{n-2}\right) \right) ,\;n%
\begin{tabular}{l}
=%
\end{tabular}%
0,1,\ldots ,  \label{eb-a} \\
q_{j}^{n} &=&-\frac{1}{W_{0}}\left( \frac{T_{j+1}^{n}-T_{j-1}^{n}}{2\Delta x}%
+\sum_{k=1}^{n}W_{k}q_{j}^{n-k}\right) ,\;n=1,2,\ldots ,  \label{ce-a}
\end{eqnarray}%
where 
\begin{equation}
W_{k}=\int_{0}^{1}\frac{\phi \left( \gamma \right) }{\left( \Delta t\right)
^{\gamma }}\omega _{k}\left( \gamma \right) \mathrm{d}\gamma
,\;\;k=0,1,\ldots ,n.  \label{w-k}
\end{equation}

Solution is found by solving the coupled system of equations (\ref{eb-a}), (%
\ref{ce-a}). In the first step ($n=0$), the initial conditions (\ref{ic-a})
used in (\ref{eb-a}) imply $T_{j}^{1}=T_{j}^{0},$ where it is assumed that $%
q_{j}^{-1}=q_{j}^{-2}=0,$ as well. The obtained $T_{j}^{1}$ is used in the
second step ($n=1$), along with the initial condition (\ref{ic-a})$_{2},$ in
(\ref{ce-a}), and $q_{j}^{1}$ is calculated. In this step, $T_{j}^{2}$ is
calculated as well, according to (\ref{eb-a}). The algorithm for marching in
time consists in repeating the calculations as in the second step. The
system of equations (\ref{eb-a}), (\ref{ce-a}) is such that, while marching
in time, the spatial domain shrinks due to the application of the leap frog
scheme in the infinite spatial domain, since the scheme does not calculate
values of solution in the outer points for each time step. The scheme
follows the rule $j=2n-1,\ldots ,J-\left( 2n-1\right) ,$ $n=1,2,\ldots ,$
for the heat flux and $j=2n,\ldots ,J-2n,$ $n=1,2,\ldots ,$ for the
temperature, where $J$ is the last point in the domain in the first step.

If, instead of using the Adams-Bashforth scheme, one uses the Euler (\ref{es}%
) and the centered scheme (\ref{cs}), the discretization of the energy
balance equation (\ref{EB-bd}), instead of (\ref{eb-a}), yields 
\begin{eqnarray}
T_{j}^{n+1} &=&T_{j}^{n}-\frac{\Delta t}{2\Delta x}\left(
q_{j+1}^{n}-q_{j-1}^{n}\right) ,\;\;n%
\begin{tabular}{l}
=%
\end{tabular}%
0,1,\ldots ,  \label{eb-a-1} \\
T_{j}^{n+1} &=&T_{j}^{n-1}-\frac{\Delta t}{\Delta x}\left(
q_{j+1}^{n}-q_{j-1}^{n}\right) ,\;\;n%
\begin{tabular}{l}
=%
\end{tabular}%
1,2,\ldots ,  \label{eb-a-2}
\end{eqnarray}%
respectively. In the case of the centered scheme with RAW filter one should
use (\ref{eb-a-2}) along with the procedure given in \cite{Williams}.
Solution is found by solving the coupled system of equations: either (\ref%
{eb-a-1}), (\ref{ce-a}) in the case of Euler discretization, or (\ref{eb-a-2}%
), (\ref{ce-a}) in the case of using the centered scheme. The coupled system
of equations (\ref{eb-a-1}), (\ref{ce-a}) is solved by following the same
procedure as described when the Adams-Bashforth scheme is used. When solving
the system of equations (\ref{eb-a-2}), (\ref{ce-a}), in the first step ($%
n=0 $), the Euler scheme is used in approximating the energy balance
equation (\ref{EB-bd}), i.e., (\ref{eb-a-1}) is assumed, implying $%
T_{j}^{1}=T_{j}^{0}, $ according to the initial conditions (\ref{ic-a}).
Following steps are same as already described for the case of using the
Adams-Bashforth scheme.

The memory effects are clearly visible in (\ref{ce-a}), since the history of
heat flux is taken into account through the weights $W_{k},$ $k=0,1,\ldots
,n.$ The relation (\ref{w-k}) shows that weights are dependent on the
constitutive model. Namely, in the case of constitutive distribution (\ref%
{CDF-bd})$_{1}$, corresponding to the multi-term time-fractional heat
conduction law, weights take the following form%
\begin{equation}
W_{k}=\frac{1}{\left( \Delta t\right) ^{\alpha _{0}}}\omega _{k}\left(
\alpha _{0}\right) +\sum_{\nu =1}^{N}\frac{\tau _{\nu }}{\left( \Delta
t\right) ^{\alpha _{\nu }}}\omega _{k}\left( \alpha _{\nu }\right)
,\;\;k=0,1,\ldots ,n,  \label{w-k-mt}
\end{equation}%
while in the case of constitutive distribution (\ref{CDF-bd})$_{2}$,
corresponding to the power-type distributed-order heat conduction law,
weights are%
\begin{equation}
W_{k}=\int_{0}^{1}\frac{1}{\left( \Delta t\right) ^{\gamma }}\omega
_{k}\left( \gamma \right) \mathrm{d}\gamma ,\;\;k=0,1,\ldots ,n.
\label{w-k-do}
\end{equation}%
It is clear that in the case of multi-term law no further approximation of
weights (\ref{w-k-mt}) is needed, as opposed to the case of power-type
distributed-order law, where the trapezoidal method for integral
calculation, used for approximating weights (\ref{w-k-do}), yields%
\begin{equation}
W_{k}=\sum_{m=0}^{M-1}\frac{1}{\left( \Delta t\right) ^{\frac{2m+1}{2}\Delta
\gamma }}\omega _{k}\left( \frac{2m+1}{2}\Delta \gamma \right) \Delta \gamma
,\;\;k=0,1,\ldots ,n,  \label{w-k-do-tr}
\end{equation}%
where $M=\frac{1}{\Delta \gamma }.$

\section{Results}

Temperature and heat flux as solutions to system of energy balance equation (%
\ref{EB-bd}) and fractional distributed-order heat conduction law (\ref%
{CE-bd}), with initial (\ref{IC-bd}) and boundary conditions (\ref{BC-bd}),
analytically obtained by (\ref{T,q}) and numerically by (\ref{eb-a}), (\ref%
{ce-a}), are plotted in cases of multi-term and power-type distributed-order
laws, represented by constitutive distribution and function (\ref{CDF-bd}).
Further, analytical response (\ref{T,q}) to the initial temperature assumed
as the Gaussian function and the solutions obtained using different
numerical schemes for energy balance equation: (\ref{eb-a}), (\ref{eb-a-1}),
and (\ref{eb-a-2}), along with the same approximation of the constitutive
equation (\ref{ce-a}), are compared.

\subsection{Analytically obtained temperature and heat flux \label{a-t-hf}}

In order to show time evolution of temperature and heat flux spatial
profiles, the initial temperature distribution is assumed as 
\begin{equation*}
T_{0}\left( x\right) =T_{0}\,\delta \left( x\right) ,\;\;x\in 
\mathbb{R}
,
\end{equation*}%
where $\delta $ is the Dirac delta distribution, so that temperature and
heat flux, given by (\ref{T,q}), reduce to the solution kernels $P$ and $Q,$
(\ref{P}) and (\ref{Q}).

In the case of multi-term heat conduction law, the model parameters are
taken as $\alpha _{0}=0,$ $\alpha _{1}=0.25,$ $\alpha _{2}=0.5,$ $\alpha
_{3}=0.75,$ $\tau _{1}=0.4,$ $\tau _{2}=0.6,$ and $\tau _{3}=0.8,$ while the
amplitude of the initial temperature distribution is $T_{0}=1.$ Figures \ref%
{fig-r-1} 
\begin{figure}[tbph]
\centering
\includegraphics[scale=0.225]{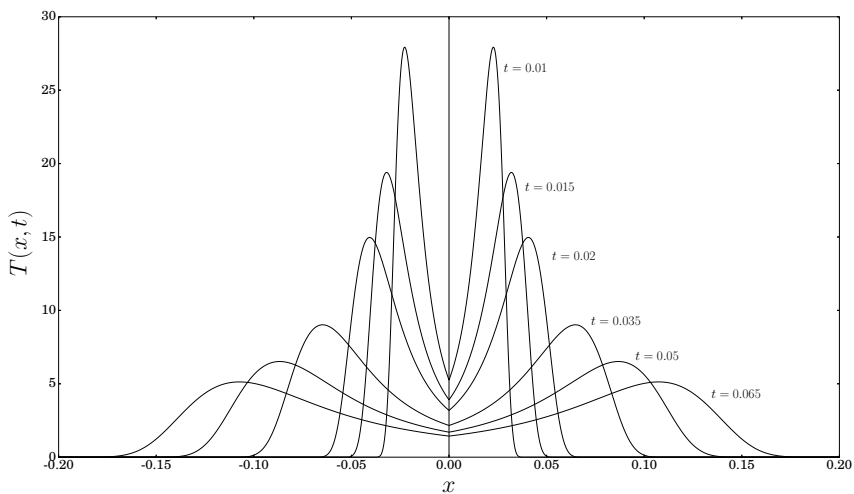}
\caption{Spatial profiles of temperature at different time instances,
obtained analitically for multi-term heat conduction law, with Dirac
distribution as initial condition. }
\label{fig-r-1}
\end{figure}
and \ref{fig-r-2} 
\begin{figure}[tbph]
\centering
\includegraphics[scale=0.225]{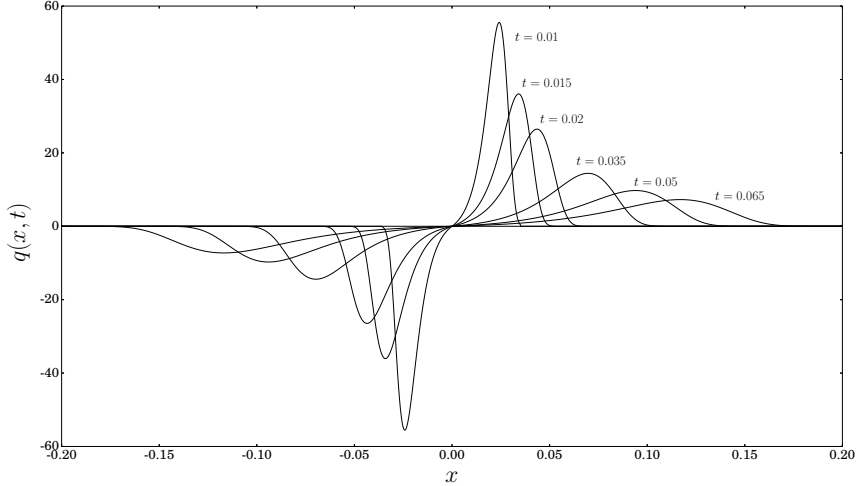}
\caption{Spatial profiles of heat flux at different time instances, obtained
analitically for multi-term heat conduction law, with Dirac distribution as
initial condition. }
\label{fig-r-2}
\end{figure}
show the spatial profiles of temperature and heat flux at different time
instances. The temperature distribution is symmetric with respect to the
vertical axis, while the heat flows away from the origin (where the initial
Dirac delta temperature distribution was introduced), therefore producing
the antisymmetric character of the heat flux spatial profiles. One notices
that spatial profiles of temperature have the similar behavior as in the
case of the wave equation with energy dissipation effects included, see for
example \cite[Figures 3.2, 3.3, 3.6, and 3.7]{APSZ-2}, unlike the heat
conduction equations with fractional Cattaneo and Jeffreys heat conduction
laws, see \cite[Figures 7.3, 7.4, and 7.14]{APSZ-1}. Namely, as time passes,
the peaks of both temperature and heat flux profiles propagate in space
(which is characteristic for wave-like behavior) and decrease in height
while increasing in width (which is characteristic for diffusion-like
behavior in which case the peaks do not propagate). Therefore, heat
conduction with multi-term heat conduction law, similarly as in the case of
the classical Cattaneo constitutive law, might be considered as the
propagation of heat waves. The diffusive characteristics of the process can
be observed through the decrease of peaks' height during time, as well as by
the increase of their width.

\begin{figure}[tbph]
\centering
\includegraphics[scale=0.225]{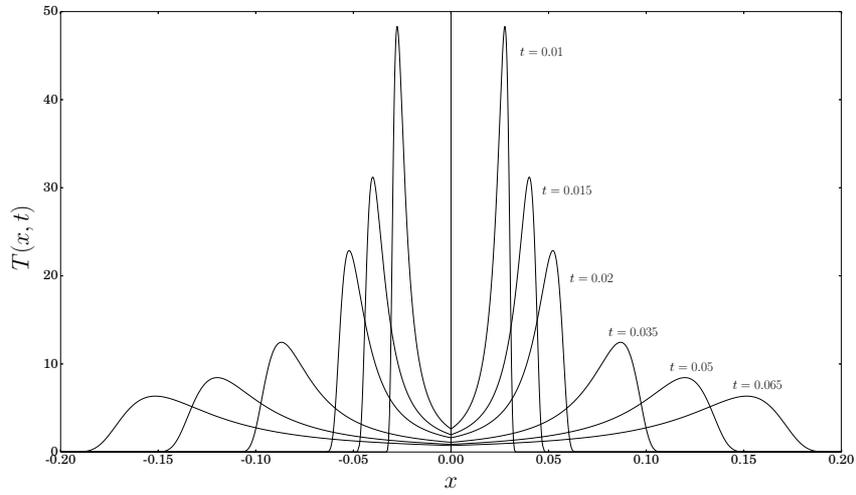}
\caption{Spatial profiles of temperature at different time instances,
obtained analitically for power-type distributed-order heat conduction law,
with Dirac distribution as initial condition.}
\label{fig-r-3}
\end{figure}

\begin{figure}[tbph]
\centering
\includegraphics[scale=0.225]{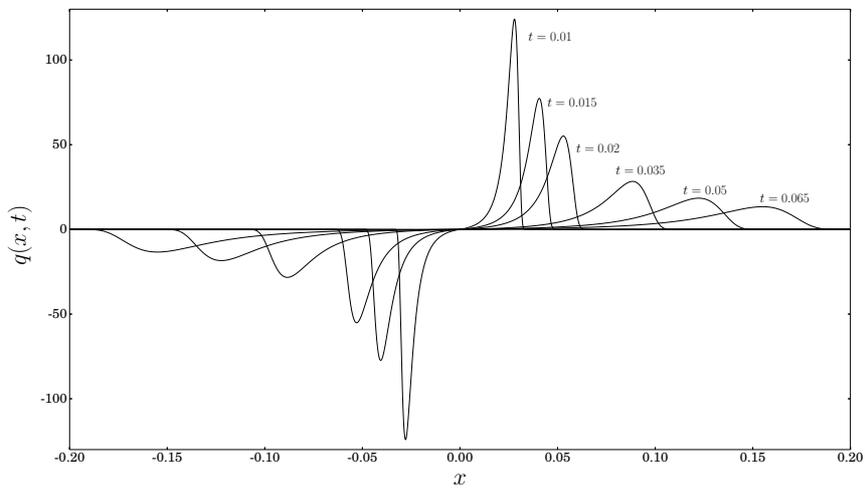}
\caption{Spatial profiles of heat flux at different time instances, obtained
analitically for power-type distributed-order heat conduction law, with
Dirac distribution as initial condition.}
\label{fig-r-4}
\end{figure}

Figures \ref{fig-r-3} 
and \ref{fig-r-4} 
show the spatial profiles of temperature and heat flux at different time
instances in the case of power-type distributed-order heat conduction law,
with $T_{0}=1.$ Likewise the case of heat conduction with multi-term law,
one might consider this type of heat conduction as the propagation of heat
waves. When compared to the case of heat conduction with multi-term law, the
wave-like character of temperature and heat flux is more prominent, since
the peaks are more localized (higher and narrower). However, the
diffusion-like character is also noticeable, since the height of the peaks
decreases while their width increases as time passes.

\subsection{Numerically obtained temperature and heat flux \label{nothf}}

In order to be able to initialize the numerical scheme for simultaneous
calculation of temperature and heat flux and to approximate the fundamental
solution to heat conduction equation, whose spatial profiles are shown in
Section \ref{a-t-hf}, the initial temperature distribution is assumed as the
Gaussian function 
\begin{equation}
T_{0}\left( x\right) =\frac{T_{0}}{2\sqrt{\pi \varepsilon }}\mathrm{e}^{-%
\frac{x^{2}}{4\varepsilon }},\;\;x\in 
\mathbb{R}
,  \label{gf}
\end{equation}%
as the smooth approximation of the Dirac delta distribution, which is
obtained as $\varepsilon \rightarrow 0.$

The aim is to compare temperature and heat flux profiles, obtained
analytically through (\ref{T,q}) by convolving the solution kernels $P$ and $%
Q,$ (\ref{P}) and (\ref{Q}), with the Gaussian function (\ref{gf}) as
initial condition, with the ones obtained numerically through (\ref{eb-a}), (%
\ref{ce-a}), with weights (\ref{w-k}) reducing to (\ref{w-k-mt}) in the case
of multi-term heat conduction law and to (\ref{w-k-do-tr}) in the case of
power-type distributed-order law.

\begin{figure}[tbph]
\centering
\includegraphics[scale=0.225]{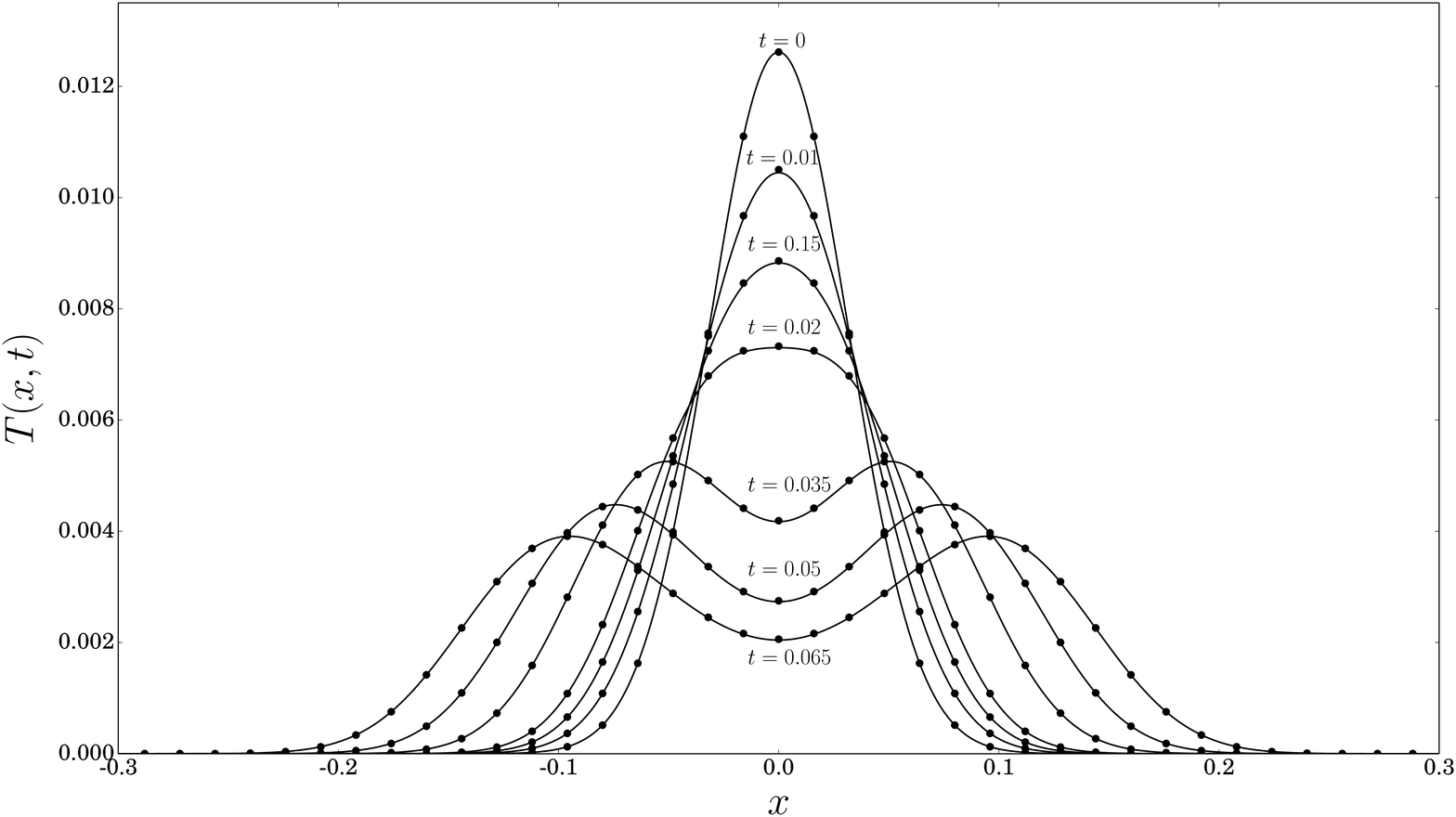}
\caption{Comparison of spatial profiles of temperature at different time
instances, obtained numerically (solid lines) and analytically (dots) for
multi-term heat conduction law, with Gaussian function as initial condition. 
}
\label{fig-r-5}
\end{figure}

\begin{figure}[tbph]
\centering
\includegraphics[scale=0.225]{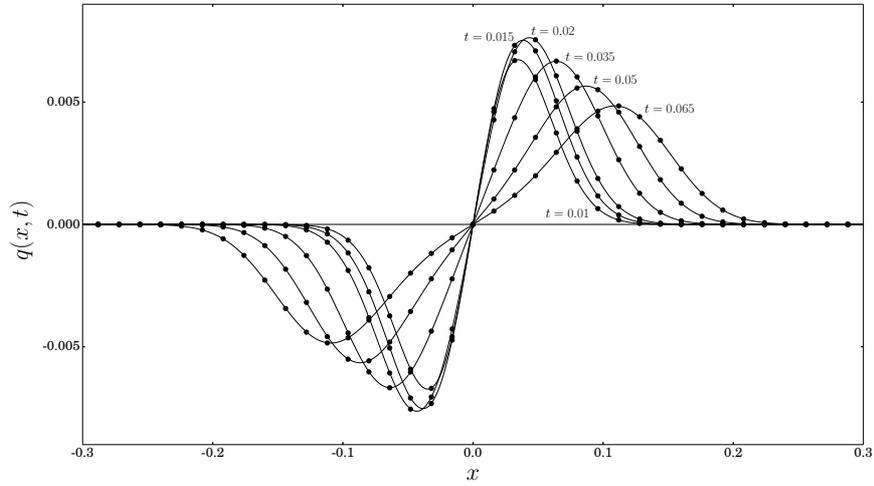}
\caption{Comparison of spatial profiles of heat flux at different time
instances, obtained numerically (solid lines) and analytically (dots) for
multi-term heat conduction law, with Gaussian function as initial condition. 
}
\label{fig-r-6}
\end{figure}

Figures \ref{fig-r-5} 
and \ref{fig-r-6} 
present the comparison of temperature and heat flux spatial profiles
obtained analytically and numerically in the case of multi-term heat
conduction law, where $T_{0}=0.001$ and $\varepsilon =0.0005$ in the
Gaussian function (\ref{gf}), the model parameters as taken as $\alpha
_{0}=0,$ $\alpha _{1}=0.25,$ $\alpha _{2}=0.5,$ $\alpha _{3}=0.75,$ $\tau
_{1}=0.4,$ $\tau _{2}=0.6,$ and $\tau _{3}=0.8,$ while the time and space
steps in the numerical scheme are $\Delta t=0.0001$ and $\Delta x=0.001.$
From Figure \ref{fig-r-5} one sees that near the initial time instant, the
temperature profiles resemble the Gaussian function and, as the process
evolves, the profiles begin to resemble the profiles of the fundamental
solution from Figure \ref{fig-r-1}, thus concluding that the initial
condition shapes the solution significantly in the beginning, while the
characteristics of the process become prominent later on. This observation
is supported by the heat flux profiles, whose peak height increases while
the Gaussian function shapes the profiles, see Figure \ref{fig-r-6}, and
decreases, as expected from fundamental solution profiles on Figure \ref%
{fig-r-2}, when the process becomes dominant.

\begin{figure}[tbph]
\centering
\includegraphics[scale=0.225]{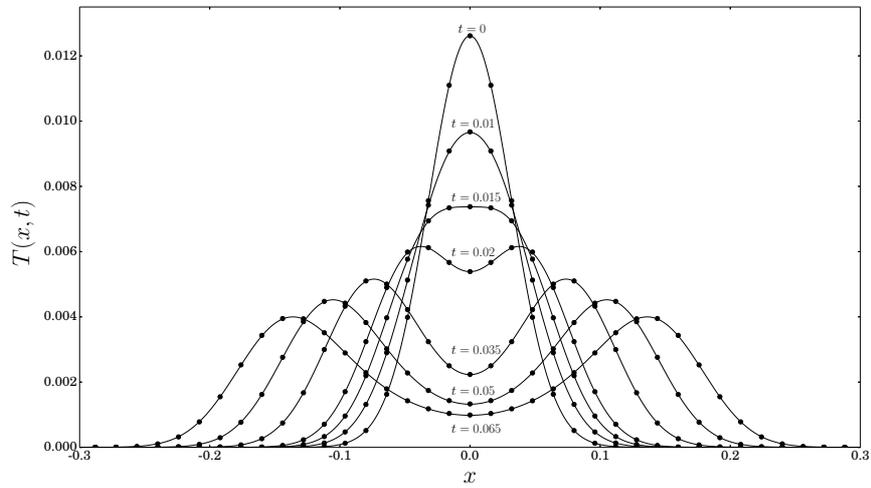}
\caption{Comparison of spatial profiles of temperature at different time
instances, obtained numerically (solid lines) and analytically (dots) for
power-type distributed-order heat conduction law, with Gaussian function as
initial condition. }
\label{fig-r-7}
\end{figure}

\begin{figure}[tbph]
\centering
\includegraphics[scale=0.225]{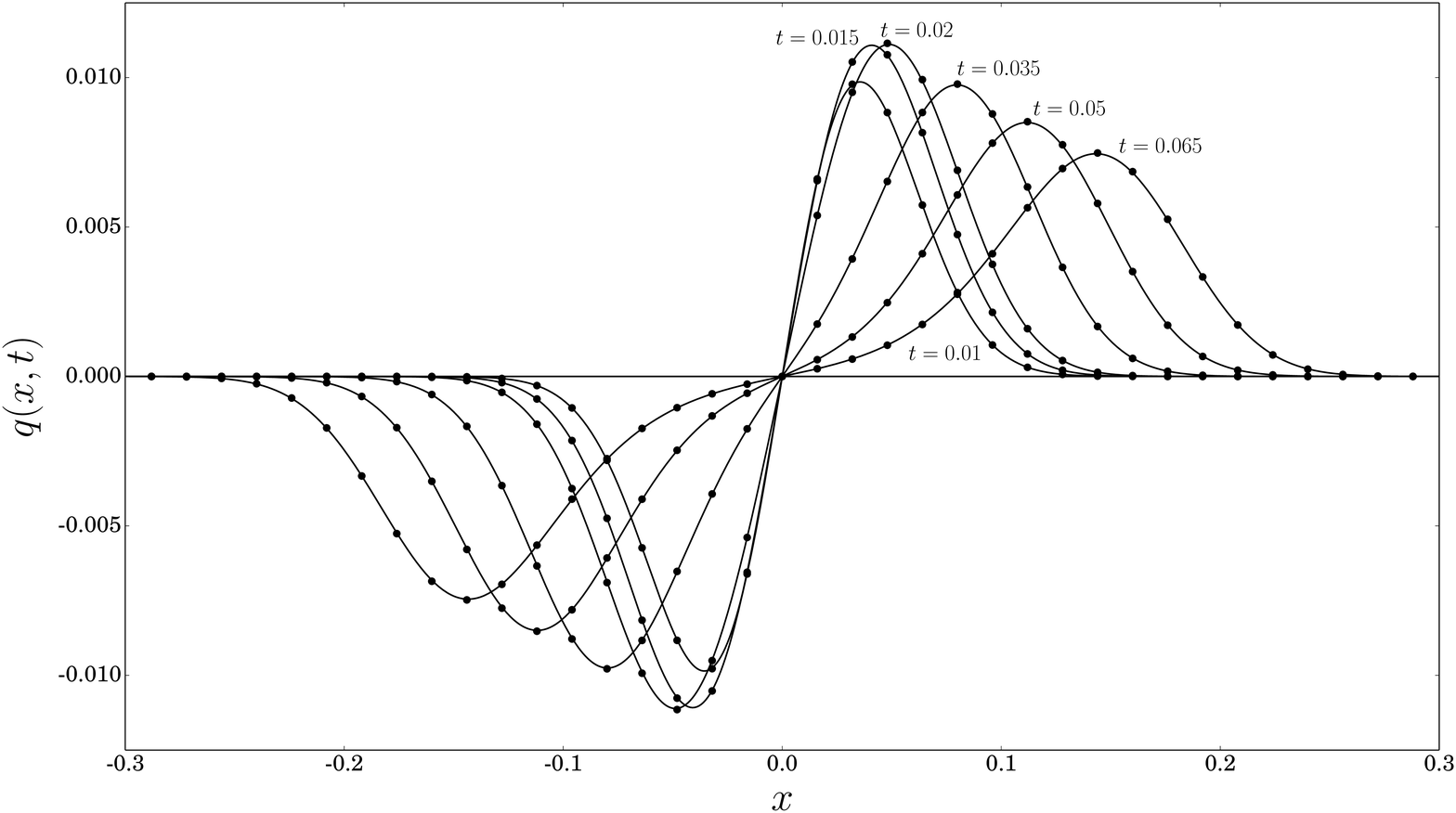}
\caption{Comparison of spatial profiles of heat flux at different time
instances, obtained numerically (solid lines) and analytically (dots) for
power-type distributed-order heat conduction law, with Gaussian function as
initial condition. }
\label{fig-r-8}
\end{figure}

In the case of power-type distributed-order law, the scheme requires an
additional discretization of the integral in weights (\ref{w-k-do}),
reducing them to (\ref{w-k-do-tr}), where the integral discretization step
is $\Delta \gamma =0.005.$ The other parameters are as in the case of
multi-term law. Comparison of temperature and heat flux spatial profiles are
shown in Figures \ref{fig-r-7} 
and \ref{fig-r-8}. 
Again, as in the case of multi-term law, near the initial time instant, the
temperature profiles resemble the Gaussian function, see Figure \ref{fig-r-7}%
, while later on, the profiles resemble the profiles of the fundamental
solution from Figure \ref{fig-r-3}. Also, the peaks' height of heat flux
profiles increase while the initial condition is dominant, see Figure \ref%
{fig-r-8}, and decrease, as fundamental one does in Figure \ref{fig-r-4},
when the process becomes dominant. As opposed to the case of multi-term law,
in this case the process begins to shape the profiles for smaller times and
the wave-like character is also more prominent.

Good agreement between analytical and numerical solution is evident from all
Figures \ref{fig-r-5} - \ref{fig-r-8}. The numerical scheme seems to be
stable for selected time length, with the discretization steps and
parameters taken.

\subsection{Comparison of numerically obtained solutions \label{cds}}

The aim is to test the accuracy and stability of numerical schemes using
different approximations of the energy balance equation (\ref{EB-bd}), while
the constitutive heat conduction law (\ref{CE-bd}) is approximated by (\ref%
{ce-a}) in all cases, by comparing the corresponding solutions among
themselves and with the analytically obtained response (\ref{T,q}) to the
Gaussian function as initial temperature. Namely, the spatial derivatives
are approximated by the leap frog scheme, the fractional derivative is used
in the Gr\"{u}nwald-Letnikov form, while the approximations of the energy
balance equation are obtained by using the following schemes:
Adams-Bashforth (\ref{eb-a}), Euler (\ref{eb-a-1}), centered (\ref{eb-a-2}),
and centered with RAW filter. The model parameters, as well as the
parameters appearing in the Gaussian function, along with the discretization
steps and time instances are throughout this section kept the same as in
Section \ref{nothf}.

\begin{figure}[tbph]
\centering
\includegraphics[scale=0.2435]{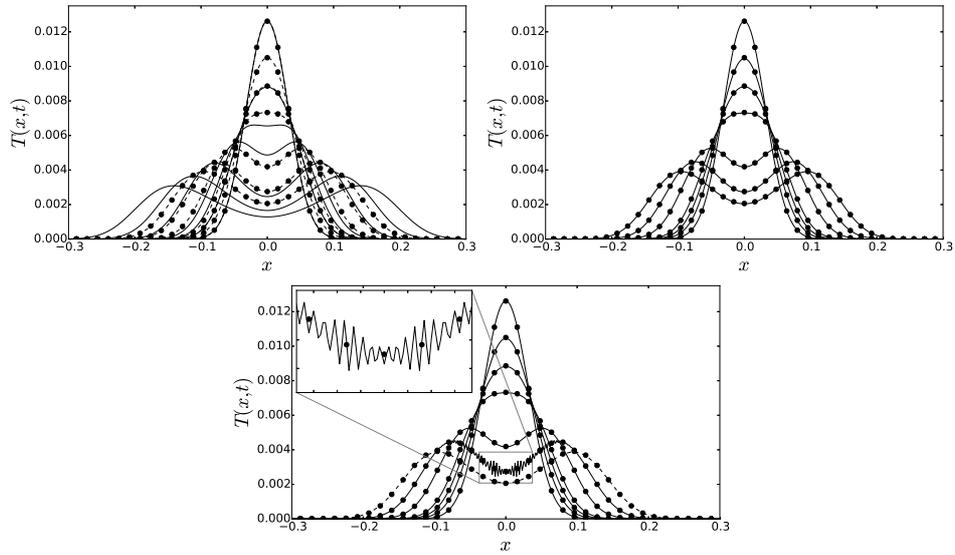}
\caption{Comparison of analytically obtained temperature profile (dots with
dashed line) with numerical one (solid line) for Euler (top left), centered
with RAW filter (top right) and centered (bottom) scheme, as response to
Gaussian function in case of multi-term heat conduction law.}
\label{fig-r-9}
\end{figure}

\begin{figure}[tbph]
\centering
\includegraphics[scale=0.2435]{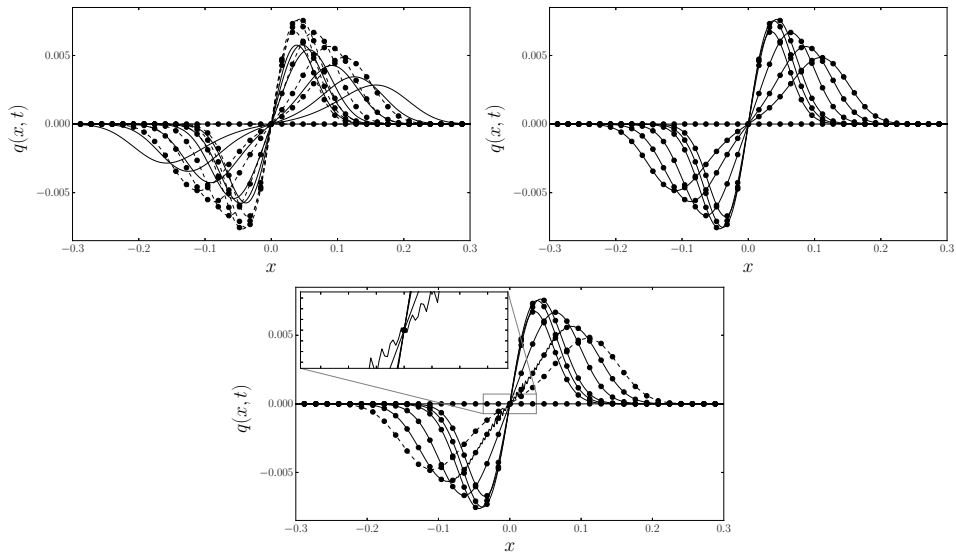}
\caption{Comparison of analytically obtained heat flux profile (dots with
dashed line) with numerical one (solid line) for Euler (top left), centered
with RAW filter (top right) and centered (bottom) scheme, as response to
Gaussian function in case of multi-term heat conduction law.}
\label{fig-r-10}
\end{figure}

Figures \ref{fig-r-9} 
and \ref{fig-r-10} 
present comparison of temperature and heat flux spatial profiles obtained
analytically according to (\ref{T,q}) and numerically according to Euler (%
\ref{eb-a-1}), centered with RAW filter, and centered (\ref{eb-a-2}) scheme
used in the energy balance equation, along with the multi-term heat
conduction law approximated by (\ref{ce-a}), (\ref{w-k-mt}), as the
responses to the Gaussian function. In the case of Euler approximation, the
numerical scheme (\ref{eb-a-1}), (\ref{ce-a}) seems to be stable, however it
yields inaccurate results for both temperature and heat flux when compared
to the analytical solution, presumably due to the first order accuracy of
the Euler scheme. Contrary to the previous case, when centered approximation
is used, numerical scheme (\ref{eb-a-2}), (\ref{ce-a}) yields accurate
results for both temperature and heat flux, however it becomes unstable by
exhibiting high frequency oscillations having small amplitudes at $t=0.05$
and having large amplitudes at $t=0.065,$ thus not depicted on Figures \ref%
{fig-r-9} and \ref{fig-r-10}. This is due to the existence of the
computational mode that may not be damped in the three level schemes, see 
\cite{Mesinger}. Having the accuracy also improved, centered scheme with RAW
filter seems to yield stable results for the parameters and time interval
considered, however it consumes more computational time when compared to
numerical scheme (\ref{eb-a-2}), (\ref{ce-a}) that uses Adams-Bashforth
approximation. Since the RAW filter averages values of the centered scheme,
it damps the amplitudes of oscillations, thus yielding the stable solution.

Tables \ref{tbl-m-t-T} 
\begin{table}[tbph]
\caption{Errors of numerically obtained solutions for temperature using
Adams-Bashforth, centered with RAW filter, centered and Euler approximation
schemes, with respect to analytically obtained response to Gaussian function
in case of multi-term heat conduction law.}
\label{tbl-m-t-T}
\begin{center}
\begin{tabular}{llllll}
& \myalign{c}{$t_{i}$} & \myalign{c}{\begin{tabular}{@{}c@{}}Adams- \\
Bashworth \end{tabular}} & \myalign{c}{\begin{tabular}{@{}c@{}} Centered
with \\ RAW filter \end{tabular}} & \myalign{c}{Centered} & %
\myalign{c}{Euler} \\ \hline
\multirow{6}{*}{$\delta_{l^{2}} \, T(t_i)$} & $0.01$ & $2.760\times 10^{-5}$
& $2.768\times 10^{-5}$ & $2.777\times 10^{-5}$ & $1.829\times 10^{-2}$ \\ 
& $0.015$ & $1.655\times 10^{-5}$ & $1.661\times 10^{-5}$ & $1.666\times
10^{-5}$ & $4.668\times 10^{-2}$ \\ 
& $0.02$ & $1.108\times 10^{-5}$ & $1.113\times 10^{-5}$ & $1.115\times
10^{-5}$ & $7.782\times 10^{-2}$ \\ 
& $0.035$ & $5.732\times 10^{-6}$ & $5.754\times 10^{-6}$ & $5.793\times
10^{-6}$ & $0.150$ \\ 
& $0.05$ & $5.167\times 10^{-6}$ & $5.174\times 10^{-6}$ & $1.020\times
10^{-3}$ & $0.203$ \\ 
& $0.065$ & $5.514\times 10^{-6}$ & $5.525\times 10^{-6}$ & $1.830\times
10^{6}$ & $0.240$ \\ \hline
\multirow{6}{*}{$\Delta_{l^{2}} \, T(t_i)$} & $0.01$ & $2.089\times 10^{-7}$
& $2.094\times 10^{-7}$ & $2.101\times 10^{-7}$ & $1.384\times 10^{-4}$ \\ 
& $0.015$ & $1.092\times 10^{-7}$ & $1.096\times 10^{-7}$ & $1.099\times
10^{-7}$ & $3.079\times 10^{-4}$ \\ 
& $0.02$ & $6.401\times 10^{-8}$ & $6.429\times 10^{-8}$ & $6.442\times
10^{-8}$ & $4.495\times 10^{-4}$ \\ 
& $0.035$ & $2.442\times 10^{-8}$ & $2.452\times 10^{-8}$ & $2.468\times
10^{-8}$ & $6.416\times 10^{-4}$ \\ 
& $0.05$ & $1.797\times 10^{-8}$ & $1.800\times 10^{-8}$ & $3.547\times
10^{-6}$ & $7.051\times 10^{-4}$ \\ 
& $0.065$ & $1.639\times 10^{-8}$ & $1.643\times 10^{-8}$ & $5.439\times
10^{3}$ & $7.144\times 10^{-4}$ \\ \hline
$\Delta_{l^{\infty}} \, T$ &  & $5.536\times 10^{-5}$ & $5.552\times 10^{-5}$
& $15.69$ & $2.456\times 10^{-3}$ \\ \hline
\end{tabular}%
\end{center}
\end{table}
and \ref{tbl-m-t-q} 
\begin{table}[tbph]
\caption{Errors of numerically obtained solutions for heat flux using
Adams-Bashforth, centered with RAW filter, centered and Euler approximation
schemes, with respect to analytically obtained response to Gaussian function
in case of multi-term heat conduction law.}
\label{tbl-m-t-q}
\begin{center}
\begin{tabular}{llllll}
& \myalign{c}{$t_{i}$} & \myalign{c}{\begin{tabular}{@{}c@{}}Adams- \\
Bashworth \end{tabular}} & \myalign{c}{\begin{tabular}{@{}c@{}} Centered
with \\ RAW filter \end{tabular}} & \myalign{c}{Centered} & %
\myalign{c}{Euler} \\ \hline
\multirow{6}{*}{$\delta_{l^{2}} \, q(t_i)$} & $0.01$ & $3.842\times 10^{-6}$
& $3.894\times 10^{-6}$ & $3.830\times 10^{-6}$ & $2.247\times 10^{-2}$ \\ 
& $0.015$ & $2.238\times 10^{-6}$ & $2.288\times 10^{-6}$ & $2.235\times
10^{-6}$ & $6.373\times 10^{-2}$ \\ 
& $0.02$ & $1.580\times 10^{-6}$ & $1.616\times 10^{-6}$ & $1.575\times
10^{-6}$ & $0.116$ \\ 
& $0.035$ & $2.521\times 10^{-6}$ & $2.563\times 10^{-6}$ & $2.507\times
10^{-6}$ & $0.254$ \\ 
& $0.05$ & $3.361\times 10^{-6}$ & $3.432\times 10^{-6}$ & $4.037\times
10^{-4}$ & $0.341$ \\ 
& $0.065$ & $3.612\times 10^{-6}$ & $3.697\times 10^{-6}$ & $7.598\times
10^{5}$ & $0.401$ \\ \hline
\multirow{6}{*}{$\Delta_{l^{2}} \, q(t_i)$} & $0.01$ & $1.458\times 10^{-8}$
& $1.478\times 10^{-8}$ & $1.454\times 10^{-8}$ & $8.528\times 10^{-5}$ \\ 
& $0.015$ & $1.152\times 10^{-8}$ & $1.178\times 10^{-8}$ & $1.151\times
10^{-8}$ & $3.282\times 10^{-4}$ \\ 
& $0.02$ & $9.128\times 10^{-9}$ & $9.334\times 10^{-9}$ & $9.101\times
10^{-9}$ & $6.708\times 10^{-4}$ \\ 
& $0.035$ & $1.370\times 10^{-8}$ & $1.393\times 10^{-8}$ & $1.363\times
10^{-8}$ & $1.378\times 10^{-3}$ \\ 
& $0.05$ & $1.491\times 10^{-8}$ & $1.522\times 10^{-8}$ & $1.790\times
10^{-6}$ & $1.514\times 10^{-3}$ \\ 
& $0.065$ & $1.306\times 10^{-8}$ & $1.336\times 10^{-8}$ & $2.746\times
10^{3}$ & $1.448\times 10^{-3}$ \\ \hline
$\Delta_{l^{\infty}} \, q$ &  & $1.426\times 10^{-5}$ & $1.437\times 10^{-5}$
& $11.19$ & $3.614\times10^{-3}$ \\ \hline
\end{tabular}%
\end{center}
\end{table}
contain absolute and relative errors of numerically obtained solutions for
temperature $T_{ns}$ and heat flux $q_{ns}$ using Adams-Bashforth (\ref{eb-a}%
), centered with RAW filter, centered (\ref{eb-a-2}) and Euler (\ref{eb-a-1}%
) approximation scheme in the energy balance equation, along with the
multi-term heat conduction law approximated by (\ref{ce-a}), (\ref{w-k-mt}),
with respect to analytically obtained responses $T_{as}$ and $q_{as}$ to the
initial Gaussian temperature distribution. The errors are calculated by
using $l^{2}$ and $l^{\infty }$ norms according to 
\begin{equation}
\Delta _{l^{2}}\,u(t_{i})=\Vert u_{as}(\cdot ,t_{i})-u_{ns}(\cdot
,t_{i})\Vert _{l^{2}},\;\;\delta _{l^{2}}\,u(t_{i})=\frac{\Delta
_{l^{2}}\,u(t_{i})}{\Vert u_{as}(\cdot ,t_{i})\Vert _{l^{2}}},\;\;\Delta
_{l^{\infty }}\,u=\Vert u_{as}-u_{ns}\Vert _{l^{\infty }},  \label{elovi}
\end{equation}%
where%
\begin{equation*}
\Vert u(\cdot ,t_{i})\Vert _{l^{2}}=\sqrt{\frac{1}{j_{\max }-j_{\min }}%
\sum_{j=j_{\min }}^{j_{\max }}\left( u(j\Delta x,t_{i})\right) ^{2}}%
,\;\;\Vert u\Vert _{l^{\infty }}=\max_{j_{\min }\leq j\leq j_{\max },\;1\leq
i\leq 6}{{u(j\Delta x,t_{i}),}}
\end{equation*}%
with $j_{\min }$ and $j_{\max }$ being dependant on the iteration step, as
described in Section \ref{fds}. For all time instances considered, the
relative and absolute $l^{2}$ as well as $l^{\infty }$ errors for both
temperature and heat flux, produced by the numerical scheme (\ref{eb-a}), (%
\ref{ce-a}), that uses the Adams-Bashforth approximation of the constitutive
law, are smaller than the corresponding errors produced by the numerical
scheme that uses the centered approximation with RAW filter. This, along
with the reduced computational time in the case of using the Adams-Bashforth
scheme implies its advantage. The scheme (\ref{eb-a-2}), (\ref{ce-a}), that
uses the centered approximation of the constitutive equation, for time
instances $t=0.01,\ldots ,0.035,$ as opposed to the temperature, produces
for the heat flux smaller values of the relative and absolute $l^{2}$ errors
than the Adams-Bashforth scheme and centered scheme with RAW filter. The
relative and absolute $l^{2}$ errors, for both temperature and heat flux,
increase when there are high frequency oscillations with small amplitudes ($%
t=0.05$) and they become significantly large when the amplitudes increase ($%
t=0.065$), which is also evident from Figures \ref{fig-r-9} and \ref%
{fig-r-10}. When using the centered scheme, the $l^{\infty }$ error is large
due to the high frequency oscillations with large amplitudes. The scheme (%
\ref{eb-a-1}), (\ref{ce-a}), that uses the Euler approximation of the
constitutive equation, produces the solution with the lowest accuracy,
except for the case when (\ref{eb-a-2}), (\ref{ce-a}) exhibits high
frequency oscillations, and the increasing error in each time instant.

\begin{figure}[tbph]
\centering
\includegraphics[scale=0.235]{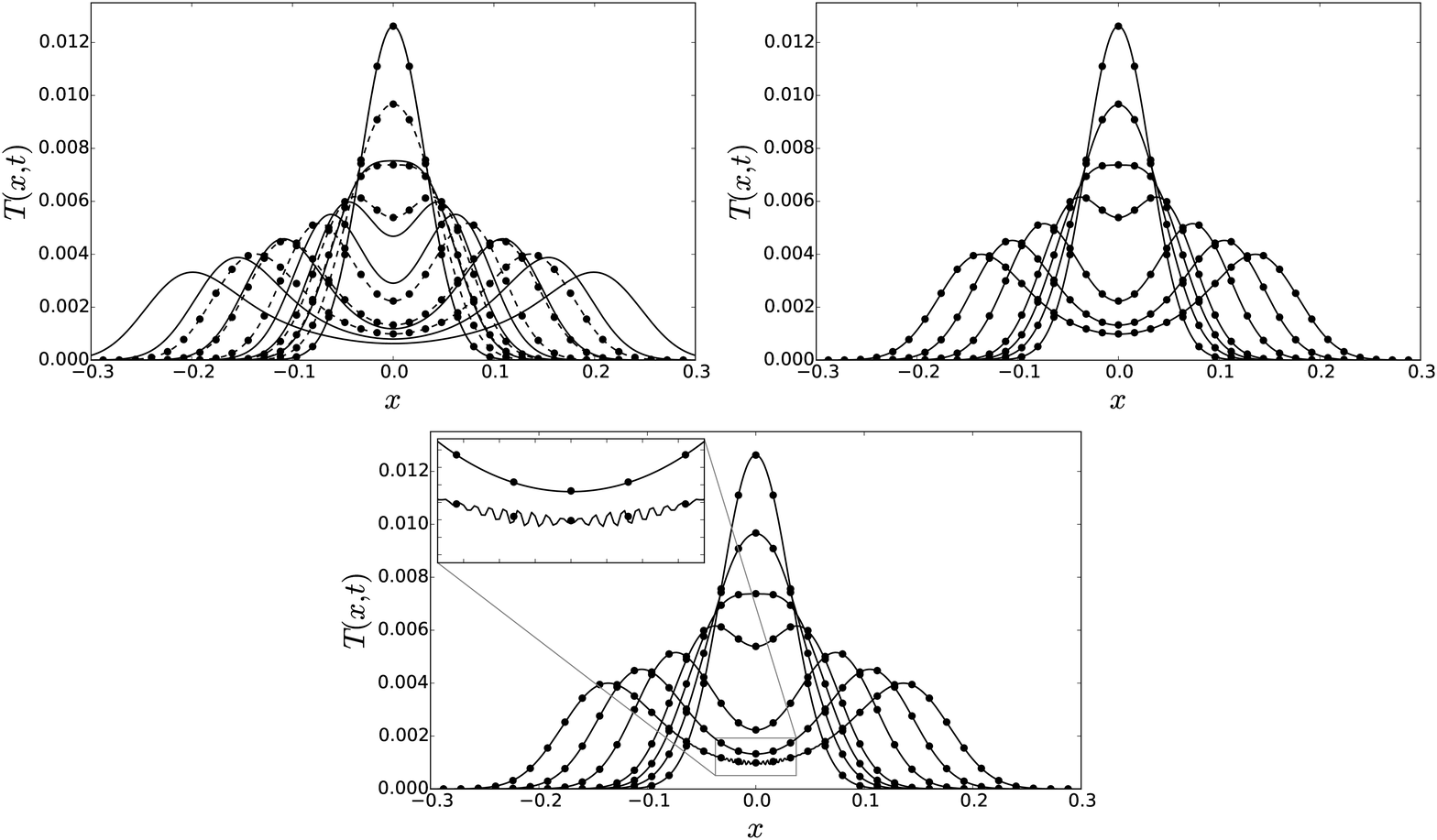}
\caption{Comparison of analytically obtained temperature profile (dots with
dashed line) with numerical one (solid line) for Euler (top left), centered
with RAW filter (top right) and centered (bottom) scheme, as response to
Gaussian function in case of power-type distributed-order heat conduction
law.}
\label{fig-r-11}
\end{figure}

\begin{figure}[tbph]
\centering
\includegraphics[scale=0.235]{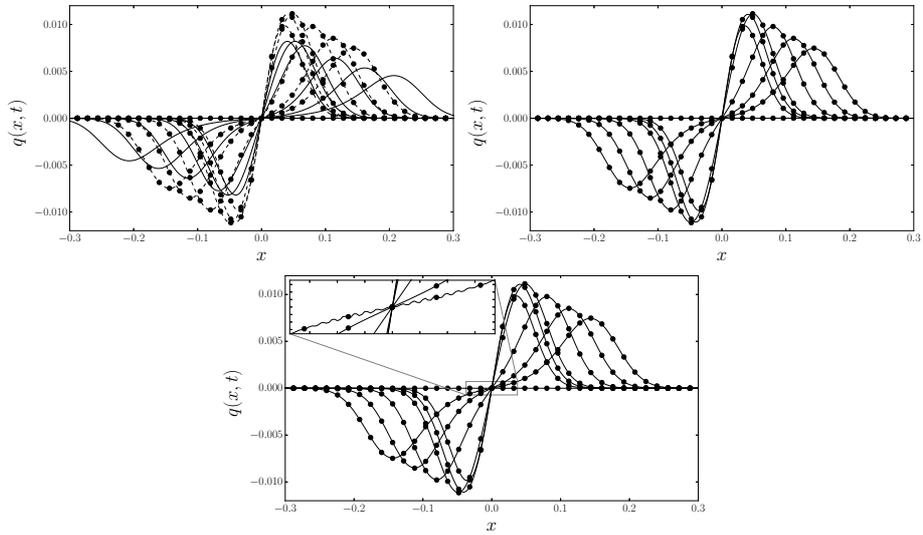}
\caption{Comparison of analytically obtained heat flux profile (dots with
dashed line) with numerical one (solid line) for Euler (top left), centered
with RAW filter (top right) and centered (bottom) scheme, as response to
Gaussian function in case of power-type distributed-order heat conduction
law.}
\label{fig-r-12}
\end{figure}

Figures \ref{fig-r-11} 
and \ref{fig-r-12} 
present comparison of temperature and heat flux spatial profiles obtained
analytically according to (\ref{T,q}) and numerically according to Euler (%
\ref{eb-a-1}), centered with RAW filter, and centered (\ref{eb-a-2}) scheme
used in the energy balance equation, along with the power-type
distributed-order heat conduction law approximated by (\ref{ce-a}), (\ref%
{w-k-do-tr}), as the responses to the Gaussian function. Similarly as in the
case of multi-term constitutive law, the use of Euler approximation gives
solutions that seem to be stable, but inaccurate; the centered scheme also
becomes unstable, but at the time instance later than the one in the case of
multi-term law; the centered scheme with RAW filter seems to yield accurate
and stable results for the parameters and time interval considered.

Tables \ref{tbl-d-o-T} 
\begin{table}[tbph]
\caption{Errors of numerically obtained solutions for temperature using
Adams-Bashforth, centered with RAW filter, centered and Euler approximation
schemes, with respect to analytically obtained response to Gaussian function
in case of power-type distributed-order heat conduction law.}
\label{tbl-d-o-T}
\begin{center}
\begin{tabular}{llllll}
& \myalign{c}{$t_{i}$} & \myalign{c}{\begin{tabular}{@{}c@{}}Adams- \\
Bashworth \end{tabular}} & \myalign{c}{\begin{tabular}{@{}c@{}} Centered
with \\ RAW filter \end{tabular}} & \myalign{c}{Centered} & %
\myalign{c}{Euler} \\ \hline
\multirow{6}{*}{$\delta_{l^{2}} \, T(t_i)$} & $0.01$ & $6.577\times 10^{-6}$
& $6.613\times 10^{-6}$ & $6.603\times 10^{-6}$ & $3.284\times 10^{-2}$ \\ 
& $0.015$ & $4.562\times 10^{-6}$ & $4.582\times 10^{-6}$ & $4.532\times
10^{-6}$ & $8.450\times 10^{-2}$ \\ 
& $0.02$ & $8.527\times 10^{-6}$ & $8.558\times 10^{-6}$ & $8.427\times
10^{-6}$ & $0.136$ \\ 
& $0.035$ & $5.579\times 10^{-6}$ & $5.643\times 10^{-6}$ & $5.542\times
10^{-6}$ & $0.273$ \\ 
& $0.05$ & $9.058\times 10^{-6}$ & $9.173\times 10^{-6}$ & $9.112\times
10^{-6}$ & $0.394$ \\ 
& $0.065$ & $1.139\times 10^{-5}$ & $1.156\times 10^{-5}$ & $7.032\times
10^{-5}$ & $0.481$ \\ \hline
\multirow{6}{*}{$\Delta_{l^{2}} \, T(t_i)$} & $0.01$ & $4.654\times 10^{-8}$
& $4.679\times 10^{-8}$ & $4.672\times 10^{-8}$ & $2.324\times 10^{-4}$ \\ 
& $0.015$ & $2.690\times 10^{-8}$ & $2.701\times 10^{-8}$ & $2.671\times
10^{-8}$ & $4.982\times 10^{-4}$ \\ 
& $0.02$ & $4.309\times 10^{-8}$ & $4.325\times 10^{-8}$ & $4.259\times
10^{-8}$ & $6.849\times 10^{-4}$ \\ 
& $0.035$ & $2.129\times 10^{-8}$ & $2.154\times 10^{-8}$ & $2.115\times
10^{-8}$ & $1.040\times 10^{-3}$ \\ 
& $0.05$ & $2.906\times 10^{-8}$ & $2.943\times 10^{-8}$ & $2.923\times
10^{-8}$ & $1.264\times 10^{-3}$ \\ 
& $0.065$ & $3.155\times 10^{-8}$ & $3.203\times 10^{-8}$ & $1.948\times
10^{-7}$ & $1.334\times 10^{-3}$ \\ \hline
$\Delta_{l^{\infty}} \, T$ &  & $2.267\times 10^{-5}$ & $2.272\times 10^{-5}$
& $8.842\times 10^{-5}$ & $2.699\times 10^{-3}$ \\ \hline
\end{tabular}%
\end{center}
\end{table}
and \ref{tbl-d-o-q} 
\begin{table}[tbph]
\caption{Errors of numerically obtained solutions for heat flux using
Adams-Bashforth, centered with RAW filter, centered and Euler approximation
schemes, with respect to analytically obtained response to Gaussian function
in case of power-type distributed-order heat conduction law.}
\label{tbl-d-o-q}
\begin{center}
\begin{tabular}{llllll}
& \myalign{c}{$t_{i}$} & \myalign{c}{\begin{tabular}{@{}c@{}}Adams- \\
Bashworth \end{tabular}} & \myalign{c}{\begin{tabular}{@{}c@{}} Centered
with \\ RAW filter \end{tabular}} & \myalign{c}{Centered} & %
\myalign{c}{Euler} \\ \hline
\multirow{6}{*}{$\delta_{l^{2}} \, q(t_i)$} & $0.01$ & $5.647\times 10^{-6}$
& $5.731\times 10^{-6}$ & $5.615\times 10^{-6}$ & $3.167\times 10^{-2}$ \\ 
& $0.015$ & $3.488\times 10^{-6}$ & $3.579\times 10^{-6}$ & $3.478\times
10^{-6}$ & $9.435\times 10^{-2}$ \\ 
& $0.02$ & $1.126\times 10^{-5}$ & $1.139\times 10^{-5}$ & $1.124\times
10^{-5}$ & $0.171$ \\ 
& $0.035$ & $8.189\times 10^{-6}$ & $8.355\times 10^{-6}$ & $8.157\times
10^{-6}$ & $0.360$ \\ 
& $0.05$ & $1.218\times 10^{-5}$ & $1.246\times 10^{-5}$ & $1.215\times
10^{-5}$ & $0.496$ \\ 
& $0.065$ & $1.481\times 10^{-5}$ & $1.517\times 10^{-5}$ & $1.989\times
10^{-5}$ & $0.595$ \\ \hline
\multirow{6}{*}{$\Delta_{l^{2}} \, q(t_i)$} & $0.01$ & $4.683\times 10^{-8}$
& $4.753\times 10^{-8}$ & $4.657\times 10^{-8}$ & $2.627\times 10^{-4}$ \\ 
& $0.015$ & $4.060\times 10^{-8}$ & $4.166\times 10^{-8}$ & $4.049\times
10^{-8}$ & $1.098\times 10^{-3}$ \\ 
& $0.02$ & $1.466\times 10^{-7}$ & $1.483\times 10^{-7}$ & $1.463\times
10^{-7}$ & $2.228\times 10^{-3}$ \\ 
& $0.035$ & $9.660\times 10^{-8}$ & $9.855\times 10^{-8}$ & $9.622\times
10^{-8}$ & $4.243\times 10^{-3}$ \\ 
& $0.05$ & $1.189\times 10^{-7}$ & $1.216\times 10^{-7}$ & $1.186\times
10^{-7}$ & $4.838\times 10^{-3}$ \\ 
& $0.065$ & $1.205\times 10^{-7}$ & $1.234\times 10^{-7}$ & $1.619\times
10^{-7}$ & $4.840\times 10^{-3}$ \\ \hline
$\Delta_{l^{\infty}} \, q$ &  & $4.375\times 10^{-5}$ & $4.399\times 10^{-5}$
& $4.368\times 10^{-5}$ & $6.074\times10^{-3}$ \\ \hline
\end{tabular}%
\end{center}
\end{table}
contain absolute and relative errors, calculated by (\ref{elovi}), of
numerically obtained solutions for temperature $T_{ns}$ and heat flux $q_{ns}
$ using Adams-Bashforth (\ref{eb-a}), centered with RAW filter, centered (%
\ref{eb-a-2}) and Euler (\ref{eb-a-1}) approximation schemes in the energy
balance equation, along with the power-type distributed-order heat
conduction law approximated by (\ref{ce-a}), (\ref{w-k-do-tr}), with respect
to analytically obtained response $T_{as}$ and $q_{as}$ to the initial
Gaussian temperature distribution. Similarly as in the case of multi-term
law, the use of Adams-Bashforth approximation produces smaller relative and
absolute $l^{2}$ and $l^{\infty }$ errors when compared with the centered
scheme with RAW filter; the Euler scheme also produces the solution with the
lowest accuracy, having relative and absolute $l^{2}$ errors increasing with
time. When compared with the Adams-Bashforth scheme and centered scheme with
RAW filter, the use of centered scheme produces smaller values of the
relative and absolute $l^{2}$ errors for both temperature and heat flux,
except at time instances $t=0.01$ and $t=0.05$ for temperature and at $%
t=0.065$ for both temperature and heat flux. As one expects, the instability
of scheme (\ref{eb-a-2}), (\ref{ce-a}) is not so clearly visible as in the
case of multi-term law, since the solution has high frequency oscillations
with small amplitudes, see Figures \ref{fig-r-11} and \ref{fig-r-12}.
However, the instability is implied by the difference of one order of
magnitude in relative $l^{2}$ error of temperature between this scheme and
the scheme that uses Adams-Bashforth approximation. This difference is not
as prominent as in the case of heat flux. In the case of temperature, the
scheme that uses Adams-Bashforth approximation has the smallest $l^{\infty }$
error, as opposed to the case of heat flux, where the smallest $l^{\infty }$
error is for the scheme that uses centered approximation.

\section{Conclusion}

The classical heat conduction equation is generalized by considering the
system of equations consisting of the energy balance equation (\ref{EB}) and
Cattaneo type time-fractional distributed-order constitutive heat conduction
law (\ref{CE}). Two cases of the constitutive equation are examined:
multi-term and power-type distributed-order heat conduction laws, with the
constitutive distribution/function given by (\ref{MT}) and (\ref{PTDO}),
respectively. The Cauchy initial value problem on the real axis is
considered by subjecting governing equations (\ref{EB}) and (\ref{CE}) to
initial and boundary conditions (\ref{IC}) and (\ref{BC}). Corresponding
dimensionless system of equations (\ref{EB-bd}) and (\ref{CE-bd}), with (\ref%
{IC-bd}) and (\ref{BC-bd}), is solved analytically through integral
transform methods: Fourier transform with respect to spatial coordinate and
Laplace transform with respect to time, as well as by the finite difference
method: leap frog numerical scheme for spatial coordinate, along with Gr\"{u}%
nwald-Letnikov and third-order Adams-Bashforth temporal numerical schemes.
The analytical solution (\ref{T,q}) is obtained as a convolution of initial
temperature distribution with the solution kernels, given by (\ref{P}) and (%
\ref{Q}), while the numerical solution is obtained through (\ref{eb-a}) and (%
\ref{ce-a}), with weights (\ref{w-k}), reducing to (\ref{w-k-mt}) and (\ref%
{w-k-do-tr}) in the cases of multi-term and power-type distributed-order
heat conduction laws, respectively. Note that solutions for temperature and
heat flux naturally arise, since the scheme requires both temperature and
heat flux for marching in time, due to the fact that the system of governing
equations is coupled.

The response to the initial Dirac delta distribution yielded temperature and
heat flux spatial profiles having the similar form as in the case of the
telegraph equation, i.e., wave equation with energy dissipation effects
included, see Figures \ref{fig-r-1} - \ref{fig-r-4}, thus describing the
propagation of heat waves, as opposed to the case of the heat conduction
equations with fractional Cattaneo and Jeffreys heat conduction laws, having
purely diffusive character.

Good agreement between analytical and numerical methods in cases of
multi-term, see Figures \ref{fig-r-5} and \ref{fig-r-6}, and power-type
distributed-order heat conduction laws, see Figures \ref{fig-r-7} and \ref%
{fig-r-8}, is found by comparing temperature and heat flux profiles,
obtained analytically by convolving the solution kernels with the Gaussian
function as initial condition and numerically through (\ref{eb-a}), (\ref%
{ce-a}), showing applicability of the joint use of Adams-Bashforth
approximation of the energy balance equation, leap frog scheme for spatial
derivatives, and Gr\"{u}nwald-Letnikov approximation of the fractional
derivative.

Justification for the use of Adams-Bashforth scheme in approximating the
energy balance equation is found by comparing the absolute and relative $%
l^{2}$ and $l^{\infty }$ errors (obtained with respect to the analytical
solutions) of temperature and heat flux, produced by using the following
schemes: Adams-Bashforth (\ref{eb-a}), Euler (\ref{eb-a-1}), centered (\ref%
{eb-a-2}), and centered with RAW filter, while the heat conduction law is in
all cases approximated by (\ref{ce-a}). Namely, the Euler scheme proved to
give stable, but inaccurate solutions, contrary to the centered scheme that
yielded unstable, but the most accurate solutions for heat flux within the
time domain of its stability, while the centered scheme with RAW filter gave
stable solutions requiring longer computational time and having higher
values of all errors. Therefore, the use of Adams-Bashforth scheme proved to
be the best choice, both because of its accuracy and stability when compared
with Euler, centered and centered with RAW filter scheme.

\appendix

\section{Justification for applicability of the Fourier inversion formula (%
\protect\ref{FIF}) \label{app}}

In order to shown that the Fourier inversion formula (\ref{FIF}) applies,
one has to prove that $s\Phi \left( s\right) \in 
\mathbb{C}
\backslash \left( -\infty ,0\right] $ for $\func{Re}s>0.$ More precisely, it
will be shown that $\arg \left( s\Phi \left( s\right) \right) \in \left[
0,\pi \right) ,$ for $\arg s\in \left[ 0,\frac{\pi }{2}\right) $ (and $\arg
\left( s\Phi \left( s\right) \right) \in \left( -\pi ,0\right) ,$ for $\arg
s\in \left( \frac{\pi }{2},0\right) $) implying that the complex square
root, $\sqrt{s\Phi \left( s\right) },$ for $\func{Re}s>0,$ is well-defined.

In the case of constitutive distribution $\phi ,$ given by (\ref{CDF-bd})$%
_{1}$, the substitution $s=\rho \,\mathrm{e}^{\mathrm{i}\varphi },$ $\rho
>0, $ $\varphi \in \left[ -\frac{\pi }{2},\frac{\pi }{2}\right] ,$ implies
that the real and imaginary parts of function $s\Phi \left( s\right) ,$ $%
\func{Re}s>0,$ where $\Phi $ is given by (\ref{cdf-lt})$_{1}$, read%
\begin{eqnarray}
\left. \func{Re}\left( s\Phi \left( s\right) \right) \right\vert _{s=\rho \,%
\mathrm{e}^{\mathrm{i}\varphi }} &=&\rho ^{\alpha _{0}+1}\cos \left( \left(
\alpha _{0}+1\right) \varphi \right) +\sum_{\nu =1}^{N}\tau _{\nu }\rho
^{\alpha _{\nu }+1}\cos \left( \left( \alpha _{\nu }+1\right) \varphi
\right) ,  \label{s-fi-re} \\
\left. \func{Im}\left( s\Phi \left( s\right) \right) \right\vert _{s=\rho \,%
\mathrm{e}^{\mathrm{i}\varphi }} &=&\rho ^{\alpha _{0}+1}\sin \left( \left(
\alpha _{0}+1\right) \varphi \right) +\sum_{\nu =1}^{N}\tau _{\nu }\rho
^{\alpha _{\nu }+1}\sin \left( \left( \alpha _{\nu }+1\right) \varphi
\right) .  \label{s-fi-im}
\end{eqnarray}%
Since, by (\ref{s-fi-im}), it holds that $\func{Im}\left( \bar{s}\Phi \left( 
\bar{s}\right) \right) =-\func{Im}\left( s\Phi \left( s\right) \right) ,$
where bar denotes the complex conjugation, it is sufficient to analyze
function $s\Phi \left( s\right) ,$ $\func{Re}s>0,$ for $\varphi \in \left[ 0,%
\frac{\pi }{2}\right] $ only. If $\varphi \in \left( 0,\frac{\pi }{2}\right]
,$ then $\func{Im}\left( s\Phi \left( s\right) \right) >0,$ since for $0\leq
\alpha _{0}<\ldots <\alpha _{N}<1,$ it is valid that implying $\sin \left(
\left( \alpha _{\nu }+1\right) \varphi \right) >0,$ $\nu =0,1,\ldots ,N.$ If 
$\varphi =0,$ then, by (\ref{s-fi-re}), $\func{Re}\left( s\Phi \left(
s\right) \right) >0.$ Therefore, $s\Phi \left( s\right) \in 
\mathbb{C}
\backslash \left( -\infty ,0\right] $ and the Fourier inversion formula (\ref%
{FIF}) applies, as well as that $\arg \left( s\Phi \left( s\right) \right)
\in \left( 0,\pi \right) ,$ for $\varphi \in \left( 0,\frac{\pi }{2}\right)
. $

\begin{figure}[tbph]
\centering
\includegraphics[scale=0.55]{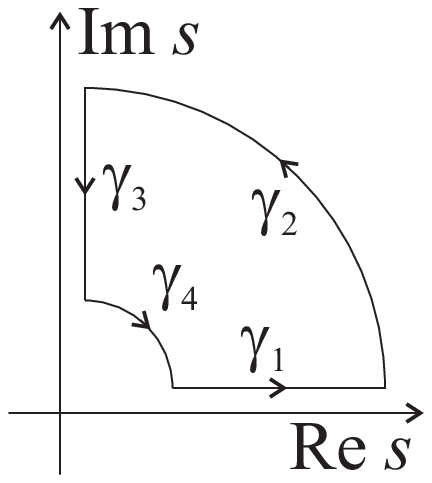}
\caption{Contour $\Gamma $.}
\label{fig-1}
\end{figure}

In the case of constitutive function $\phi ,$ given by (\ref{CDF-bd})$_{2},$
it will be shown using the argument principle that function 
\begin{equation}
\psi \left( s\right) =s\Phi \left( s\right) +\xi ^{2},\;\;s\in 
\mathbb{C}
,  \label{psi}
\end{equation}%
where $\Phi $ is given by (\ref{cdf-lt})$_{2}$, has no zeros in the right
complex half-plane ($\func{Re}s>0$) for any $\xi \in 
\mathbb{R}
,$ implying the applicability of the Fourier inversion formula (\ref{FIF}).
Moreover, it will also be shown that $\arg \left( \psi \left( s\right)
\right) \in \left( 0,\pi \right) ,$ for $\varphi \in \left( 0,\frac{\pi }{2}%
\right) .$ By substituting $s=\rho \,\mathrm{e}^{\mathrm{i}\varphi },$ $\rho
>0,$ $\varphi \in \left[ -\frac{\pi }{2},\frac{\pi }{2}\right] ,$ in (\ref%
{psi}), the real and imaginary parts of function $\psi $ are obtained as%
\begin{eqnarray}
\func{Re}\psi \left( \rho ,\varphi \right) &=&\rho \frac{\ln \rho \left(
\rho \cos \left( 2\varphi \right) -\cos \varphi \right) +\varphi \left( \rho
\sin \left( 2\varphi \right) -\sin \varphi \right) }{\ln ^{2}\rho +\varphi
^{2}}+\xi ^{2},  \label{psi-re} \\
\func{Im}\psi \left( \rho ,\varphi \right) &=&\rho \frac{\ln \rho \left(
\rho \sin \left( 2\varphi \right) -\sin \varphi \right) -\varphi \left( \rho
\cos \left( 2\varphi \right) -\cos \varphi \right) }{\ln ^{2}\rho +\varphi
^{2}}.  \label{psi-im}
\end{eqnarray}%
Similarly as in the previous case, $\func{Im}\left( \psi \left( \bar{s}%
\right) \right) =-\func{Im}\left( \psi \left( s\right) \right) ,$ thus is
sufficient to analyze function $\psi $ only in the upper right complex
quarter-plane. In order to apply the argument principle, the contour $\Gamma
=\gamma _{1}\cup \gamma _{2}\cup \gamma _{3}\cup \gamma _{4},$ shown in
Figure \ref{fig-1}, is used. 
The contour $\gamma _{1}$ is parameterized by $s=x$, $x\in (r,R),$ with $%
r\rightarrow 0$ and $R\rightarrow \infty $, so that function $\psi ,$ (\ref%
{psi}), reads%
\begin{equation*}
\psi \left( x\right) =\frac{x\left( x-1\right) }{\ln x}+\xi ^{2}>0,
\end{equation*}%
since $x-1$ and $\ln x$ are of the same sign for $x\in \left( 0,\infty
\right) .$ Moreover, 
\begin{equation*}
\psi (x)\rightarrow \xi ^{2}\;\;\text{as}\;\;x\rightarrow 0\;\;\text{and}%
\;\;\psi (x)\rightarrow \infty \;\;\text{as}\;\;x\rightarrow \infty .
\end{equation*}%
The contour $\gamma _{2}$ is parametrized by $s=R\,\mathrm{e}^{\mathrm{i}%
\varphi }$, $\varphi \in \left( 0,\frac{\pi }{2}\right) ,$ with $%
R\rightarrow \infty $. For $R$ sufficiently large, by (\ref{psi-re}) and (%
\ref{psi-im}), it is obtained%
\begin{equation*}
\func{Re}\psi \left( R,\varphi \right) \sim \frac{R^{2}}{\ln R}\cos \left(
2\varphi \right) +\xi ^{2}\;\;\text{and}\;\;\func{Im}\psi \left( R,\varphi
\right) \sim \frac{R^{2}}{\ln R}\sin \left( 2\varphi \right) >0.
\end{equation*}%
In particular, 
\begin{gather*}
\func{Re}\psi \left( R,0\right) \rightarrow \infty \;\;\text{and}\;\;\func{Re%
}\psi \left( R,\frac{\pi }{2}\right) \rightarrow -\infty ,\;\;\text{as}%
\;\;R\rightarrow \infty , \\
\func{Im}\psi \left( R,0\right) =0\;\;\text{and}\;\;\func{Im}\psi \left( R,%
\frac{\pi }{2}\right) \rightarrow \infty ,\;\;\text{as}\;\;R\rightarrow
\infty .
\end{gather*}%
Along $\gamma _{3}$, which is parametrized by $s=\rho \,\mathrm{e}^{\mathrm{i%
}\frac{\pi }{2}}$, $\rho \in \left( r,R\right) $, with $r\rightarrow 0$ and $%
R\rightarrow \infty $, by (\ref{psi-re}) and (\ref{psi-im}), the real and
imaginary parts of $\psi $ read%
\begin{equation*}
\func{Re}\psi \left( \rho ,\frac{\pi }{2}\right) =-\rho \frac{\rho \ln \rho +%
\frac{\pi }{2}}{\ln ^{2}\rho +\frac{\pi ^{2}}{4}}+\xi ^{2}\;\;\text{and}\;\;%
\func{Im}\psi \left( \rho ,\frac{\pi }{2}\right) =\rho \frac{\frac{\pi }{2}%
\rho -\ln \rho }{\ln ^{2}\rho +\frac{\pi ^{2}}{4}}>0,
\end{equation*}%
since $\rho >\ln \rho ,$ for all $\rho \in \left( 0,\infty \right) .$ Also,%
\begin{gather*}
\func{Re}\psi \left( \rho ,\frac{\pi }{2}\right) \rightarrow \xi ^{2},\;\;%
\text{as}\;\;\rho \rightarrow 0\;\;\text{and}\;\;\func{Re}\psi \left( \rho ,%
\frac{\pi }{2}\right) \rightarrow -\infty ,\;\;\text{as}\;\;\rho \rightarrow
\infty , \\
\func{Im}\psi \left( \rho ,\frac{\pi }{2}\right) \rightarrow 0,\;\;\text{as}%
\;\;\rho \rightarrow 0\;\;\text{and}\;\;\func{Im}\psi \left( \rho ,\frac{\pi 
}{2}\right) \rightarrow \infty ,\;\;\text{as}\;\;\rho \rightarrow \infty .
\end{gather*}%
The last part of the contour $\Gamma $ is the arc $\gamma _{4}$,
parametrized by $s=r\,\mathrm{e}^{\mathrm{i}\varphi }$, $\varphi \in \left(
0,\frac{\pi }{2}\right) $, with $r\rightarrow 0$. Again, (\ref{psi-re}) and (%
\ref{psi-im}), for $r$ sufficiently small, yield%
\begin{equation*}
\func{Re}\psi \left( r,\varphi \right) =-\frac{r}{\ln r}\cos \varphi +\xi
^{2}>0\;\;\text{and}\;\;\func{Im}\psi \left( r,\varphi \right) =-\frac{r}{%
\ln r}\sin \varphi >0,
\end{equation*}%
as well as%
\begin{gather*}
\func{Re}\psi \left( r,0\right) \rightarrow \xi ^{2}\;\;\text{and}\;\;\func{%
Re}\psi \left( r,\frac{\pi }{2}\right) \rightarrow \xi ^{2},\;\;\text{as}%
\;\;r\rightarrow 0, \\
\func{Im}\psi \left( r,0\right) =0\;\;\text{and}\;\;\func{Im}\psi \left( r,%
\frac{\pi }{2}\right) \rightarrow 0,\;\;\text{as}\;\;r\rightarrow 0.
\end{gather*}

Summing up, it is evident that as the complex variable $s$ changes along the
contour $\Gamma ,$ with $r$ tending to zero and $R$ tending to infinity, the
imaginary part of function $\psi ,$ (\ref{psi}), stays non-negative implying
that $\arg \left( \psi \left( s\right) \right) \in \left( 0,\pi \right) ,$
for $\varphi \in \left( 0,\frac{\pi }{2}\right) .$ This ensures the
applicability of the Fourier inversion formula (\ref{FIF}) and that the
complex square root of $\psi $ is well-defined.

\section*{Acknowledgment}

This work is supported by the Serbian Ministry of Science, Education and
Technological Development under grant $174005$, as well as by the Provincial
Government of Vojvodina under grant $114-451-2098$.

Figures presenting plots of solutions have been produced using Matplotlib, 
\cite{Hunter:2007}.


\begin{thebibliography}{10}

\bibitem{Alvarez}
J.~Alvarez-Ramirez, G.~Fernandez-Anaya, F.~J. Valdes-Parada, and J.~A.
  Ochoa-Tapia.
\newblock A high-order extension for the {C}attaneo's diffusion equation.
\newblock {\em Physica A}, 368:345--354, 2006.

\bibitem{AKOZ-1}
T.~M. Atanackovic, S.~Konjik, Lj. Oparnica, and D.~Zorica.
\newblock The {C}attaneo type space-time fractional heat conduction equation.
\newblock {\em Continuum Mechanics and Thermodynamics}, 24:293--311, 2012.

\bibitem{APSZ-1}
T.~M. Atanackovic, S.~Pilipovic, B.~Stankovic, and D.~Zorica.
\newblock {\em Fractional Calculus with Applications in Mechanics: Vibrations
  and Diffusion Processes}.
\newblock Wiley-ISTE, London, 2014.

\bibitem{APSZ-2}
T.~M. Atanackovic, S.~Pilipovic, B.~Stankovic, and D.~Zorica.
\newblock {\em Fractional Calculus with Applications in Mechanics: Wave
  Propagation, Impact and Variational Principles}.
\newblock Wiley-ISTE, London, 2014.

\bibitem{APZ}
T.~M. Atanackovic, S.~Pilipovic, and D.~Zorica.
\newblock Diffusion wave equation with two fractional derivatives of different
  order.
\newblock {\em Journal of Physics A: Mathematical and Theoretical},
  40:5319--5333, 2007.

\bibitem{APZ-1}
T.~M. Atanackovic, S.~Pilipovic, and D.~Zorica.
\newblock Time distributed-order diffusion-wave equation. {I}. {V}olterra type
  equation.
\newblock {\em Proceedings of the Royal Society A: Mathematical, Physical and
  Engineering Sciences}, 465:1869--1891, 2009.

\bibitem{APZ-2}
T.~M. Atanackovic, S.~Pilipovic, and D.~Zorica.
\newblock Time distributed-order diffusion-wave equation. {II}. {A}pplications
  of the {L}aplace and {F}ourier transformations.
\newblock {\em Proceedings of the Royal Society A: Mathematical, Physical and
  Engineering Sciences}, 465:1893--1917, 2009.

\bibitem{BorinoDiPaolaZingales}
G.~Borino, M.~{Di Paola}, and M.~Zingales.
\newblock A non-local model of fractional heat conduction in rigid bodies.
\newblock {\em European Physical Journal Special Topics}, 193:173--184, 2011.

\bibitem{CamargoOliveiraVaz}
R.~F. Camargo, E.~{Capelas de Oliveira}, and J.~{Vaz Jr.}
\newblock On the generalized {M}ittag-{L}effler function and its application in
  a fractional telegraph equation.
\newblock {\em Mathematical Physics Analysis and Geometry}, 15:1--16, 2012.

\bibitem{CamargoCO08}
R.~F. Camargo, A.~O. Chiacchio, and E.~Capelas de~Oliveira.
\newblock Differentiation to fractional orders and the fractional telegraph
  equation.
\newblock {\em Journal of Mathematical Physics}, 49:033505--1--12, 2008.

\bibitem{CascavalEcksteinFrotaGoldstein}
R.~C. Cascaval, E.~C. Eckstein, C.~L. Frota, and J.~A. Goldstein.
\newblock Fractional telegraph equations.
\newblock {\em Journal of Mathematical Analysis and Applications},
  276:145--159, 2002.

\bibitem{ChenLiuAnh}
J.~Chen, F.~Liu, and V.~Anh.
\newblock Analytical solution for the time-fractional telegraph equation by the
  method of separating variables.
\newblock {\em Journal of Mathematical Analysis and Applications},
  338:1364--1377, 2008.

\bibitem{CompteMetzler}
A.~Compte and R.~Metzler.
\newblock The generalized {C}attaneo equation for the description of anomalous
  transport processes.
\newblock {\em Journal of Physics A: Mathematical and General}, 30:7277--7289,
  1997.

\bibitem{Durran}
D.~R. Durran.
\newblock The third-order {A}dams-{B}ashforth method: {A}n attractive
  alternative to leapfrog time differencing.
\newblock {\em Monthly Weather Review}, 119:702--720, 1991.

\bibitem{FernandezAnaya}
G.~Fernandez-Anaya, F.~J. Valdes-Parada, and J.~Alvarez-Ramirez.
\newblock On generalized fractional {C}attaneo's equations.
\newblock {\em Physica A}, 390:4198--4202, 2011.

\bibitem{GLM}
R.~Gorenflo, Y.~Luchko, and F.~Mainardi.
\newblock Wright functions as scale-invariant solutions of the diffusion-wave
  equation.
\newblock {\em Journal of Computational and Applied Mathematics}, 118:175--191,
  2000.

\bibitem{Hanyga-sf}
A.~Hanyga.
\newblock Multidimensional solutions of space-fractional diffusion equations.
\newblock {\em Proceedings of the Royal Society A: Mathematical, Physical and
  Engineering Sciences}, 457:2993--3005, 2001.

\bibitem{Hanyga-stf}
A.~Hanyga.
\newblock Multi-dimensional solutions of space-time-fractional diffusion
  equations.
\newblock {\em Proceedings of the Royal Society A: Mathematical, Physical and
  Engineering Sciences}, 458:429--450, 2002.

\bibitem{Hanyga-tf}
A.~Hanyga.
\newblock Multidimensional solutions of time-fractional diffusion-wave
  equations.
\newblock {\em Proceedings of the Royal Society A: Mathematical, Physical and
  Engineering Sciences}, 458:933--957, 2002.

\bibitem{HuLiuAnhTurner}
X.~Hu, F.~Liu, V.~Anh, and I.~Turner.
\newblock A numerical investigation of the time distributed-order diffusion
  model.
\newblock {\em ANZIAM}, 55 (EMAC2013):C464--C478, 2014.

\bibitem{HuLiuTurnerAnh}
X.~Hu, F.~Liu, I.~Turner, and V.~Anh.
\newblock An implicit numerical method of a new time distributed-order and
  two-sided space-fractional advection-dispersion equation.
\newblock {\em Numerical Algorithms}, 72:393--407, 2016.

\bibitem{Huang}
F.~Huang.
\newblock Analytical solution for the time-fractional telegraph equation.
\newblock {\em Journal of Applied Mathematics}, 2009:890158--1--9, 2009.

\bibitem{Hunter:2007}
J.~D. Hunter.
\newblock Matplotlib: {A} 2{D} graphics environment.
\newblock {\em Computing in Science \& Engineering}, 9:90--95, 2007.

\bibitem{Jos}
D.~D. Joseph and L.~Preziosi.
\newblock Heat waves.
\newblock {\em Reviews of Modern Physics}, 61:41--73, 1989.

\bibitem{TAFDE}
A.~A. Kilbas, H.~M. Srivastava, and J.~J. Trujillo.
\newblock {\em Theory and Applications of Fractional Differential Equations}.
\newblock Elsevier B.V., Amsterdam, 2006.

\bibitem{LiuZhengLiuZhang}
L.~Liu, L.~Zheng, F.~Liu, and X.~Zhang.
\newblock An improved heat conduction model with {R}iesz fractional
  {C}attaneo-{C}hristov flux.
\newblock {\em International Journal of Heat and Mass Transfer},
  103:1191--1197, 2016.

\bibitem{Mai96}
F.~Mainardi.
\newblock Fractional relaxation-oscillation and fractional diffusion-wave
  phenomena.
\newblock {\em Chaos, Solitons \& Fractals}, 7:1461--1477, 1996.

\bibitem{MaLuPa}
F.~Mainardi, Y.~Luchko, and G.~Pagnini.
\newblock The fundamental solution of the space-time fractional diffusion
  equation.
\newblock {\em Fractional Calculus and Applied Analysis}, 4:153--192, 2001.

\bibitem{MaGoMi}
F.~Mainardi, A.~Mura, R.~Gorenflo, and M.~Stojanovic.
\newblock The two forms of fractional relaxation of distributed order.
\newblock {\em Journal of Vibration and Control}, 13:1249--1268, 2007.

\bibitem{MainardiPagniniGorenflo}
F.~Mainardi, G.~Pagnini, and R.~Gorenflo.
\newblock Some aspects of fractional diffusion equations of single and
  distributed order.
\newblock {\em Applied Mathematics and Computation}, 187:295--305, 2007.

\bibitem{MaPaMuGo}
F.~Mainardi, G.~Pagnini, A.~Mura, and R.~Gorenflo.
\newblock Time-fractional diffusion of distributed order.
\newblock {\em Journal of Vibration and Control}, 14:1267--1290, 2008.

\bibitem{MeerschaertNaneVellaisamy}
M.~M. Meerschaert, E.~Nane, and P.~Vellaisamy.
\newblock Distributed-order fractional diffusions on bounded domains.
\newblock {\em Journal of Mathematical Analysis and Applications},
  379:216--228, 2011.

\bibitem{Mesinger}
F.~Mesinger and A.~Arakawa.
\newblock {\em Numerical methods used in atmospheric models}, volume~1 of {\em
  GARP publications series No. 17}.
\newblock World Meteorological Organization, International Council of
  Scientific Unions, Geneve, 1976.

\bibitem{MishraRai}
T.~N. Mishra and K.~N. Rai.
\newblock Numerical solution of {FSPL} heat conduction equation for analysis of
  thermal propagation.
\newblock {\em Applied Mathematics and Computation}, 273:1006--1017, 2016.

\bibitem{Mohanty}
R.~K. Mohanty.
\newblock An unconditionally stable difference scheme for the
  one-space-dimensional linear hyperbolic equation.
\newblock {\em Applied Mathematics Letters}, 17:101--105, 2004.

\bibitem{MongioviZingales}
M.~S. Mongiov\'{i} and M.~Zingales.
\newblock A non-local model of thermal energy transport: {T}he fractional
  temperature equation.
\newblock {\em International Journal of Heat and Mass Transfer}, 67:593--601,
  2013.

\bibitem{MorgadoRebelo}
M.~L. Morgado and M.~Rebelo.
\newblock Numerical approximation of distributed order reaction-diffusion
  equations.
\newblock {\em Journal of Computational and Applied Mathematics}, 275:216--227,
  2015.

\bibitem{pod}
I.~Podlubny.
\newblock {\em Fractional Differential Equations}.
\newblock Academic Press, San Diego, 1999.

\bibitem{QiGuo}
H.~Qi and X.~Guo.
\newblock Transient fractional heat conduction with generalized {C}attaneo
  model.
\newblock {\em International Journal of Heat and Mass Transfer}, 76:535--539,
  2014.

\bibitem{QiJiang}
H.~Qi and X.~Jiang.
\newblock Solutions of the space-time fractional {C}attaneo diffusion equation.
\newblock {\em Physica A}, 390:1876--1883, 2011.

\bibitem{QiXuGuo}
HT. Qi, HY. Xu, and XW. Guo.
\newblock The {C}attaneo-type time fractional heat conduction equation for
  laser heating.
\newblock {\em Computers and Mathematics with Applications}, 66:824--831, 2013.

\bibitem{RJ}
M.~R. Rapai\'{c} and Z.~D. Jeli\v{c}i\'{c}.
\newblock Optimal control of a class of fractional heat diffusion systems.
\newblock {\em Nonlinear Dynamics}, 62:39--51, 2010.

\bibitem{ShenLiuAnh}
S.~Shen, F.~Liu, and V.~Anh.
\newblock Numerical approximations and solution techniques for the space-time
  {R}iesz-{C}aputo fractional advection-diffusion equation.
\newblock {\em Numerical Algorithms}, 56:383--403, 2011.

\bibitem{ShenLiuLiuAnh}
S.~Shen, F.~Liu, Q.~Liu, and V.~Anh.
\newblock Numerical simulation of anomalous infiltration in porous media.
\newblock {\em Numerical Algorithms}, 68:443--454, 2015.

\bibitem{SunWu}
ZZ. Sun and X.~Wu.
\newblock A fully discrete difference scheme for a diffusion-wave system.
\newblock {\em Applied Numerical Mathematics}, 56:193--209, 2006.

\bibitem{ZecovaTerpak}
M.~\v{Z}ecov\'{a} and J.~Terp\'{a}k.
\newblock Heat conduction modeling by using fractional-order derivatives.
\newblock {\em Applied Mathematics and Computation}, 257:365--373, 2015.

\bibitem{Williams}
P.~D. Williams.
\newblock A proposed modification to the {R}obert-{A}sselin time filter.
\newblock {\em Monthly Weather Review}, 137:2538--2546, 2009.

\bibitem{XuQiJiang}
HY. Xu, HT. Qi, and XY. Jiang.
\newblock Fractional {C}attaneo heat equation in a semi-infinite medium.
\newblock {\em Chinese Physics B}, 22:014401--1--6, 2013.

\bibitem{YakubovichRodrigues}
S.~Yakubovich and M.~M. Rodrigues.
\newblock Fundamental solutions of the fractional two-parameter telegraph
  equation.
\newblock {\em Integral Transforms and Special Functions}, 23:509--519, 2012.

\bibitem{YeLiuAnh}
H.~Ye, F.~Liu, and V.~Anh.
\newblock Compact difference scheme for distributed-order time-fractional
  diffusion-wave equation on bounded domains.
\newblock {\em Journal of Computational Physics}, 298:652--660, 2015.

\bibitem{YeLiuAnhTurnerAMC}
H.~Ye, F.~Liu, V.~Anh, and I.~Turner.
\newblock Maximum principle and numerical method for the multi-term time-space
  {R}iesz-{C}aputo fractional differential equations.
\newblock {\em Applied Mathematics and Computation}, 227:531--540, 2014.

\bibitem{YeLiuAnhTurnerIMA}
H.~Ye, F.~Liu, V.~Anh, and I.~Turner.
\newblock Numerical analysis for the time distributed-order and {R}iesz space
  fractional diffusions on bounded domains.
\newblock {\em IMA Journal of Applied Mathematics}, 80:825--838, 2015.

\bibitem{YusteQuintanaMurillo}
S.~B. Yuste and J.~{Quintana-Murillo}.
\newblock Fast, accurate and robust adaptive finite difference methods for
  fractional diffusion equations.
\newblock {\em Numerical Algorithms}, 71:207--228, 2016.

\bibitem{ZhuangLiu}
P.~Zhuang and F.~Liu.
\newblock Implicit difference approximation for the time fractional diffusion
  equation.
\newblock {\em Journal of Applied Mathematics and Computing}, 22:87--99, 2006.

\bibitem{Zingales}
M.~Zingales.
\newblock Fractional-order theory of heat transport in rigid bodies.
\newblock {\em Communications in Nonlinear Science and Numerical Simulation},
  19:3938--3953, 2014.

\end{thebibliography}

\end{document}